\documentclass[letterpaper,12pt,openright,twoside]{amsart}
\pdfoutput=1

\usepackage{ifthen}
\newboolean{PrintVersion}
\setboolean{PrintVersion}{false}

\usepackage[foot]{amsaddr}
\usepackage[refpage]{nomencl} 
\usepackage{amsmath,amssymb,amstext} 
\usepackage[pdftex]{graphicx} 
\usepackage{amsthm}
\usepackage{thm-restate}
\usepackage{aliascnt}
\usepackage{enumerate}
\usepackage{bbm}
\usepackage{svg}
\usepackage[pdftex,letterpaper=true,pagebackref=false]{hyperref} 
\hypersetup{
    plainpages=false,       
    pdfpagelabels=true,     
    bookmarks=true,         
    unicode=false,          
    pdftoolbar=true,        
    pdfmenubar=true,        
    pdffitwindow=false,     
    pdfstartview={FitH},    
    pdftitle={Unavoidable minors and tree-decompositions},    
    pdfauthor={Jim Geelen and Benson Joeris},    
    pdfsubject={Graph Theory},  
    pdfkeywords={graph theory} {connectivity} {minor} {tree-width} {tree-decomposition}, 
    pdfnewwindow=true,      
    colorlinks=true,        
    linkcolor=blue,         
    citecolor=green,        
    filecolor=magenta,      
    urlcolor=cyan           
}
\ifthenelse{\boolean{PrintVersion}}{   
\hypersetup{	
    citecolor=black,%
    filecolor=black,%
    linkcolor=black,%
    urlcolor=black}
}{} 

\newcommand{\tryincludegraphics}[2][]{%
\IfFileExists{#2.eps}{%
  \includegraphics[#1]{#2}%
}{%
 \IfFileExists{#2.pdf}{%
  \includegraphics[#1]{#2}%
  }{ 
   \IfFileExists{#2.jpg}{%
    \includegraphics[#1]{#2}%
    }{
      \IfFileExists{#2.png}{%
      \includegraphics[#1]{#2}%
      }{%
      \begin{tcolorbox}[width=6cm,height=4cm,arc=0mm,auto outer arc]
      \end{tcolorbox}
     }
   }
 }
}%
}

\makeatletter
\def\newaliasedtheorem#1[#2]#3{%
  \newaliascnt{#1@alt}{#2}
  \newtheorem{#1}[#1@alt]{#3}
  \expandafter\newcommand\csname #1@altname\endcsname{#3}
}
\makeatother
\newtheorem{theorem}{Theorem}[section]
\newaliasedtheorem{lemma}[theorem]{Lemma}
\newaliasedtheorem{corollary}[theorem]{Corollary}
\newaliasedtheorem{proposition}[theorem]{Proposition}
\newaliasedtheorem{conjecture}[theorem]{Conjecture}
\newtheorem{claim}{Claim}[theorem]
\theoremstyle{definition}

\newtheorem*{question*}{Question}
\newtheorem*{theorem*}{Theorem}

\newtheorem*{case-one}{Case 1}
\newtheorem*{case-two}{Case 2}
\newtheorem*{case-three}{Case 3}

\newenvironment{subproof}{\proof[Proof of Claim]}{\endproof}

\newcommand{\N}{\mathbb{N}}
\newcommand{\Z}{\mathbb{Z}}

\newcommand{\cA}{\mathcal{A}}
\newcommand{\cB}{\mathcal{B}}

\newcommand{\cM}{\mathcal{M}}

\newcommand{\cP}{\mathcal{P}}
\newcommand{\cQ}{\mathcal{Q}}

\newcommand{\cX}{\mathcal{X}}
\newcommand{\cY}{\mathcal{Y}}

\renewcommand{\phi}{\varphi}
\DeclareMathOperator{\tw}{tw}

\newcommand{\delete}{\,\backslash\,}

\newcommand{\fnum}[1]{\hyperref[#1]{\ensuremath{f_{\ref*{#1}}}}}
\makeatletter
\newcommand{\fcur}{\ensuremath{f_{\@currentlabel}}}
\makeatother

\title{A Generalization of the Grid Theorem}
\author{Jim Geelen \and Benson Joeris}
\address{Department of Combinatorics and Optimization,
University of Waterloo, Waterloo, Canada} 
\thanks{This research was partially supported by a grant from the
Office of Naval Research [N00014-10-1-0851].}

\subjclass[2010]{05C40 (Primary), 05C83 (Secondary)}
\keywords{connectivity, tree-width, unavoidable minors, grid theorem}
\date{\today}

\begin{document}
\begin{abstract}
  A graph has tree-width at most $k$ if it can be obtained from a set of graphs
  each with at most $k+1$ vertices by a sequence of clique sums. We refine this
  definition by, for each non-negative integer $\theta$, defining the
  $\theta$-tree-width of a graph to be at most $k$ if it can be obtained from a
  set of graphs each with at most $k+1$ vertices by a sequence of clique sums on
  cliques of size less than $\theta$. We find the unavoidable minors for the
  graphs with large $\theta$-tree-width and we obtain Robertson and Seymour's
  Grid Theorem as a corollary.
\end{abstract}

\maketitle

\section{Introduction}
We introduce a refinement of tree-width and prove a generalization of Robertson
and Seymour's Grid Theorem; see~\cite{GM-V}. We start with the statement of our
main result which gives an unavoidable-minor characterization of graphs with
large ``$\theta$-connected sets'' and then we relate this to tree-width.

For a non-negative integer $\theta$, a set $Z$ of vertices in a graph $G$ is
called {\em $\theta$-connected} if for each pair $(X,Y)$ of subsets of $Z$ with
$|X|=|Y|\le \theta$, $G$ contains $|X|$ vertex-disjoint paths between $X$ and
$Y$; some of these paths may be edgeless, when $X\cap Y\neq\emptyset$.

\subsection*{Generalized wheels}
\begin{figure}
  \centering
  \includesvg[clean,width=0.8\columnwidth]{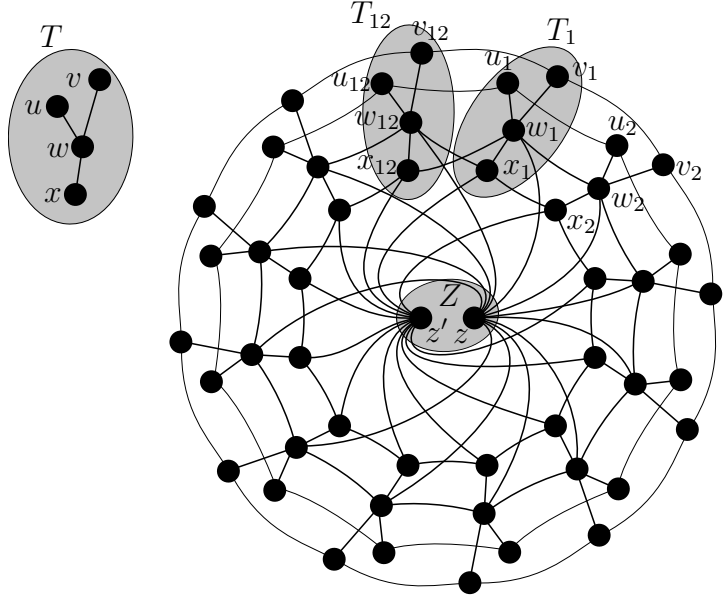}
  \caption{A $(4,2,12)$-wheel together with its tree $T$. Note that the matching
    between $T_{12}$ and $T_1$ differs from the matchings between the other
    copies of $T$. Also note that $z$ is adjacent to each copy of $w$ while $z'$
    is adjacent to each copy of $x$.}
  \label{fig:wheel-defn}
\end{figure}

Next we describe a class of graphs which have large $\theta$-connected sets;
see~\autoref{fig:wheel-defn}.

Let $t,\ell,n\in\N$ with $t\geq 1$ and $n\geq 3$.
Now let $T$ be a $t$-vertex tree, let $Z$ be an $\ell$-element
set, let $\pi:V(T)\rightarrow V(T)$ be a 
permutation, and let $\psi: Z\rightarrow V(T)$ be an arbitrary function.
Then the \emph{$(t,\ell,n)$-wheel defined by $(T,Z,\pi,\psi)$}
is the graph $G$ constructed as follows:
\begin{enumerate}[(W$_1$)]
\item Let $G_1$ be the disjoint union of $n$ copies of $T$, named
  $T_1,\ldots,T_n$. For each $v\in V(T)$ and $i\in\{1,\ldots,n\}$, let $v_i$
  denote the copy of $v$ in $T_i$.
\item Let $G_2$ be obtained from $G_1$ by adding an edge between $v_i$ and
  $v_{i+1}$ for each $v\in V(T)$ and each $i\in\{1,\ldots,n-1\}$.
\item Let $G_3$ be obtained from $G_2$ by adding an edge between $v_n$ and
  $(\pi(v))_1$ for each $v\in V(T)$.
\item Let $G$ be obtained from $G_3$ by adding $Z$ as a set of isolated vertices
  and then, for each $z\in Z$ and $i\in\{1,\ldots,n\}$, adding an edge from $z$
  to $(\psi(z))_i$.
\end{enumerate}

A {\em $(\theta;n)$-wheel} is a $(t,\ell,n)$-wheel where $2t+\ell = \theta$.
In~\autoref{thm:wheel-conn-set} we prove that each $(\theta;n)$-wheel contains a
$\theta$-connected set with $n$ vertices. Also note that if $n\ge \theta$ then
the complete bipartite graph $K_{\theta,n}$ contains a $\theta$-connected set
with $n$ vertices.

Our main result is the following.
\begin{theorem}
  \label{thm:main-exists}
  For each $\theta,n\in\N$ with $\theta\ge 2$ and $n\geq 3$, there exists
  $m\in\N$ such that, if $G$ is a graph with $\theta$-tree-width at least $m$,
  then $G$ contains either a $K_{\theta,n}$-minor or a $(\theta;n)$-wheel-minor.
\end{theorem}

In fact, we prove a slight strengthening of \autoref{thm:main-exists},
\autoref{thm:main-local}. This is obtain as a corollary of
\autoref{thm:nearly-balanced-necklace}, which is interesting in its own
right\textemdash it shows that $\theta$-connected sets gives rise to
\emph{subgraphs} in a family closely related to the $(\theta;n)$-wheels. These
results appeared in the doctoral dissertation of the second author;
see~\cite{joeris-thesis}.

Oporowski, Oxley, and Thomas~\cite{Oporowski1993} found the ``unavoidable
minors'' for large $3$-connected graphs and for large internally $4$-connected
graphs; both of those results can be deduced as corollaries
of~\autoref{thm:main-exists}.

\subsection*{Tree-width}
We will give the more usual definition of tree-width in terms of tree
decompositions in~\autoref{sec:k-connected-set-duality}; for brevity we define
tree-width in terms of ``clique-sums'' here. Let $G_1$ and $G_2$ be graphs with
a common complete subgraph $H$. Suppose that $V(G_1)\cap V(G_2) = V(H)$ and that
$E(G_1)\cap E(G_2) = E(H)$. Then we refer to the graph $(G_1\cup G_2) - E(H)$ as
the {\em clique-sum} of $G_1$ and $G_2$; the {\em order} of the sum is $|V(H)|$.

The {\em tree-width} of a graph $G$ is the smallest integer $k$ such that $G$
can be obtained from a set of graphs each with at most $k+1$ vertices by a
sequence of clique-sums. We introduce the following refinement: for a
non-negative integer $\theta$, the {\em $\theta$-tree-width} of $G$ is the
smallest integer $k$ such that $G$ can be obtained from a set of graphs each
with at most $k+1$ vertices by a sequence of clique-sums of order less than
$\theta$.

There is a well-known connection between large tree-width in graphs and the
existence of large highly-connected vertex-sets; this connection was first
observed by Robertson, Seymour and Thomas~\cite{Robertson1994a} and later
refined by Diestel, Jensen, Gorbonov and Thomassen~\cite{Diestel1999} who proved
that: for each $\omega\in\N$ and each graph $G$,
\begin{enumerate}[(i)]
\item if $G$ contains an $(\omega+1)$-connected set of size at least $3\omega$,
  then $G$ has tree-width at least $\omega$, and
\item conversely, if $G$ has no $(\omega+1)$-connected set of size at least
  $3\omega$, then $G$ has tree-width less than $4\omega$.
\end{enumerate}

We prove the following result which relates large $\theta$-tree-width in a graph
to the existence of a large $\theta$-connected set.

\begin{theorem}
  \label{thm:k-tree-width-connected-set-duality}
  For each integer $\theta\geq 3$, if $U$ is a maximum cardinality
  $\theta$-connected set in a graph $G$, and $\omega$ is the $\theta$-tree-width
  of $G$, then
  \[\omega < |U| \leq \binom{\omega+1}{\theta-1}(\theta-1).\]
\end{theorem}

In~\autoref{sec:gridtheorem} we derive Robertson and Seymour's Grid Theorem
from~\autoref{thm:main-exists}; our proof of \autoref{thm:main-exists} does not
make use of the Grid Theorem, and hence provides an alternative proof of the
Grid Theorem.

\section{Preliminaries}
\label{separations}

In this section we redefine $\theta$-connected set, tree-width, and
$\theta$-tree-width in term of ``separations''; these definitions are easily
shown to be equivalent to the definitions in the introduction, but are more
convenient for our proofs. We also prove that each $(\theta;n)$-wheel contains a
$\theta$-connected set with $n$ vertices.

\subsection*{Separations}

A \emph{separation} in a graph $G$ is an ordered pair of subgraphs $(G_1,G_2)$
such that $G= G_1\cup G_2$; the \emph{order} of a separation $(G_1,G_2)$,
denoted $\lambda(G_1,G_2)$, is $|V(G_1)\cap V(G_2)|$. For $\theta\in\N$, a
$\theta$-separation is a separation of order at most $\theta$. The following
submodular inequality is both well-known and easy to prove: if $(G_1,G_2)$ and
$(H_1,H_2)$ are separations of $G$, then $(G_1\cap H_1,G_2\cup H_2)$ and
$(G_1\cup H_1,G_2\cap H_2)$ are both separations of $G$ and
$$ \lambda(G_1\cap H_1,G_2\cup H_2)
+ \lambda(G_1\cup H_1,G_2\cap H_2) \le 
\lambda(G_1,G_2) + \lambda(H_1,H_2).
$$

For an integer $\theta\ge 0$, a set $Z$ of vertices in a graph $G$ is called
{\em $\theta$-connected} if
$\min\{|U\cap V(G_1)|,|U\cap V(G_2)|\}\leq\lambda(G_1,G_2)$ for each separation
$(G_1,G_2)$ of order less than $\theta$ in $G$.

\subsection*{Paths}

An \emph{end} of a path $P$ is a vertex of degree at most one in $P$. An
\emph{$(X,Y)$-path} is a path whose set of ends is $\{x,y\}$, with $x\in X$ and
$y\in Y$; note that $x$ and $y$ need not be distinct.
Two paths $P$ and $P'$ in a graph $G$ are \emph{internally-disjoint} if each vertex in $P\cap P'$
is an end of both paths.

\subsection*{Rim-transversals}
Consider a $(t,\ell,n)$-wheel $H$ defined by the tuple $(T,Z,\pi,\psi)$. The
tree $T$ is called the \emph{rim-tree}, the sequence $(T_1,\ldots,T_n)$ is
called the \emph{rim-tree-sequence}, and the vertices in $Z$ are called the
\emph{hubs} of $H$. A {\em rim-transversal} of $H$ is a set $U\subseteq V(H)-Z$
that contains exactly one vertex from each of the trees $T_1,\ldots,T_n$.

\begin{theorem}
  \label{thm:wheel-conn-set}
  For $t,\ell,n\in\N$, if $U$ is a rim-transversal of a $(t,\ell,n)$-wheel $H$,
  then $U$ is a $(2t+\ell)$-connected set in $H$.
  \begin{proof}
    Let $(T_1,\ldots,T_n)$ be the rim-tree-sequence of $H$. Consider any
    $(2t+\ell-1)$-separation $(H_1,H_2)$ in $H$. Let $X=V(H_1\cap H_2)$, let $I$
    be the set of integers $i\in\{1,\ldots,n\}$ such that
    $X\cap V(T_i) = \emptyset$.

    Consider any $i,j\in I$. By definition, $X$ is disjoint from $T_i$ and
    $T_j$, and, by the definition of the $(t,\ell,n)$-wheel $H$, there exist a
    family $\cP$ of $2t+\ell$ distinct, pairwise internally-disjoint
    $(V(T_i),V(T_j))$-paths. Then each $x\in X$ is in at most one path of $\cP$,
    and $|X|<|\cP|$, so the vertex set $V(T_i\cup T_j)$ is contained in a single
    connected component of $H-X$. Hence, up to symmetry, we may assume that
    $V(T_i)\subseteq V(H_1)$ for each $i\in I$.

    For each $u\in U\cap V(H_2)$, $u\in V(T_i)$ for some $i\in\{1,\ldots,n\}-I$,
    so $V(T_i)\cap X\neq\emptyset$. Therefore,
    $|V(H_2)\cap U|\leq |X|= \lambda(H_1,H_2).$ It follows that $U$ is a
    $(2t+\ell)$-connected set in $H$.
  \end{proof}
\end{theorem}

\subsection*{Tree decompositions}
A {\em leaf} of a tree $T$ is a degree-one vertex.
A \emph{tree-decomposition} of a graph $G$ is a tree $T$ such that
$E(G)$ is a subset of the leaves of $T$;
we say that an edge $t$ of $G$ {\em labels} the leaf $t$ of $T$. We will refer
to the vertices of a tree-decomposition as \emph{nodes}, to avoid confusion with
the vertices of $G$. For a vertex $v$ of $G$ we let $\delta_G(v)$ denote the set
of edges of $G$ that are incident with $v$ and let $T_v$ denote the unique
minimal subtree of $T$ containing those leaves of $T$ that are labelled by the
edges in $\delta_G(v)$. For a node $t\in V(T)$ we define a {\em bag}
$Y_t\subseteq V(G)$ to be the set of all vertices $v\in V(G)$ such that
$t\in V(T_v)$. The \emph{width} of a tree-decomposition $T$ is defined as
$\max\{|Y_t|-1\, :\, t\in V(T)\}$. The \emph{tree-width} of a graph $G$ is the
minimum width over all tree-decompositions of $G$.

The {\em bag} associated with an edge $e=tt'$ of a tree-decomposition $T$ of a
graph $G$, denoted $Y_e$, is defined as $Y_t\cap Y_{t'}$ and the {\em adhesion}
of $e$ is defined as $|Y_e|$. Let $T_1$ and $T_2$ be the components of $T-e$
and, for each $i\in\{1,2\}$, let $G_i$ be the subgraph of $G$ induced by the
edges of $G$ labeling $T_i$. Then $(G_1,G_2)$ is a separation of order $|Y_e|$;
we refer to this as the separation of $G$ {\em associated} with $e$. The {\em
  adhesion} of $T$ is defined as $\max\{|Y_e|\,:\,e\in E(T)\}$.

For each natural number $\theta$, the \emph{$\theta$-tree-width} of a graph $G$,
denoted $\tw_\theta(G)$, is the minimum width over all tree-decompositions
of $G$ of adhesion less than $\theta$.

When there is more than one tree-decomposition under consideration we use
$Y_e(T)$ and $Y_v(T)$ to denote the bags $Y_e$ and $Y_v$. We will use the terms
{\em node-bags} and {\em edge-bags} to distinguish the two types of bags.

\section{Tree-width and highly-connected sets}
\label{sec:k-connected-set-duality}
In this section we prove \autoref{thm:k-tree-width-connected-set-duality} which
shows that a graph $G$ has large $\theta$-tree-width if and only if $G$ has a
large $\theta$-connected set.
\begin{theorem}
  \label{thm:conn-set-decomposition}
  For each natural number $\theta$, each graph $G$ contains a tree-decomposition
  $T$ of adhesion less than $\theta$ and width at most $\tw_\theta(G)$ such that
  each node-bag is a $\theta$-connected set in $G$.
  \begin{proof}
    Let $\omega=\tw_\theta(G)$. For a tree-decomposition $T$ of $G$ and for each
    $k\in\{0,\ldots,\omega+1\}$ we let $n_k(T)$ denote the number of $k$-element
    node-bags in $T$. Choose a tree-decomposition $T$ of $G$ of adhesion less
    than $\theta$ and width equal to $\omega$ such that
    $(n_{\omega+1}(T),n_{\omega}(T),\ldots,n_0(T))$ is lexicographically
    minimum. Suppose for contradicting that there is a node of $T$ whose bag is
    not $\theta$-connected, and among all such nodes we choose $r\in V(T)$
    maximizing $|Y_{r}|$.

    For each node $t\in V(T)-\{r\}$, let $e_t$ be the unique edge incident with
    $t$ such that $t$ and $r$ are in distinct components of $T\delete e_t$, let
    $T_t$ be the component of $T\delete e_t$ that contains $t$, let $A_t$ be the
    subgraph of $G$ induced by the edges labelling $T_t$, and let $B_t$ be the
    subgraph of $G$ induced by the edges labelling $T-V(T_t)$. Note that
    $(A_t,B_t)$ is a separation in $G$ of order $|Y_{e_t}|\le \theta-1$; see
    \autoref{fig:conn-set-decomposition-proof}.

    \begin{figure}
      \centering
      \includesvg[clean,width=0.3\columnwidth]{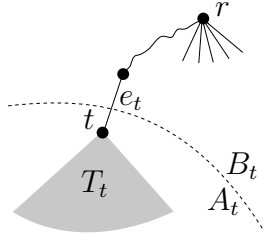}
      \caption{Notation used in the proof of
        \autoref{thm:conn-set-decomposition}.}
      \label{fig:conn-set-decomposition-proof}
    \end{figure}

    A subgraph $H\subseteq G$ is said to \emph{nest} in a separation $(A,B)$ of
    $G$ if either $V(H)\subseteq V(A)$ or $V(H)\subseteq V(B)$.

    Recall that the bag $Y_{r}$ is not $\theta$-connected. Now among all
    $(\theta-1)$-separations $(G_1,G_2)$ with
    $|V(G_1)\cap Y_r|>\lambda(G_1,G_2)$ and $|V(G_2)\cap Y_r|>\lambda(G_1,G_2)$
    we choose one minimizing $\lambda(G_1,G_2)$ and, subject to that, we choose
    one minimizing the number of nodes $t\in V(T)-\{r\}$ for which $A_t$ does
    not nest in $(G_1,G_2)$. We treat the indices of $(G_1,G_2)$ as elements of
    the cyclic group $\Z_2$; thus $G_{i+2} = G_i$.

    \begin{claim}
      \label{thm:conn-set-decomposition-1}
      For each $t\in V(T)-\{r\}$ and $i\in\Z_2$,
      \[\lambda(A_t\cap G_i,B_t \cup G_{i+1}) \leq \lambda(A_t,B_t),\]
      and if equality holds then $A_t$ nests in $(G_1,G_2)$.
      \begin{subproof}
        Let $t\in V(T)-\{r\}$. Without loss of generality, we may assume $i=1$.
        Suppose that
        \[\lambda(A_t\cap G_1,B_t \cup G_2) \geq \lambda(A_t,B_t).\]
        By submodularity,
        \[\lambda(A_t\cup G_1,B_t\cap G_2) \leq \lambda(G_1,G_2).\]
        Also,
        $$|Y_{r}\cap V(A_t \cup G_1)|\ge |Y_{r}\cap V(G_1)|>\lambda(G_1,G_2)$$
        and, since $Y_{r}\subseteq V(B_t)$,
        $$|Y_{r}\cap V(B_t\cap G_2)|=|Y_{r}\cap V(G_2)|>\lambda(G_1,G_2).$$
        By the choice of $(G_1,G_2)$,
        \[\lambda(A_t\cup G_1,B_t \cap G_2) = \lambda(G_1,G_2)\]
        and, by submodularity,
        \[\lambda(A_t\cap G_1,B_t \cup G_2) = \lambda(A_t,B_t).\]

        Suppose for contradiction that $A_t$ does not nest in $(G_1,G_2)$. Note
        that $A_t$ nests in $(A_t\cup G_1,B_t\cap G_2)$ so, by the choice of
        $(G_1,G_2)$, there exists some $t'\in V(T)-\{r\}$ such that $A_{t'}$
        nests in $(G_1,G_2)$ but not in $(A_t\cup G_1,B_t\cap G_2)$. Note that
        $V(A_t)\not\subseteq V(A_{t'})$, or else $A_t$ would nest in $(G_1,G_2)$ as
        $A_{t'}$ does, and $V(A_{t'})\not\subseteq V(A_t)$, or else
        $V(A_{t'})\subseteq V(A_t\cup G_1)$. Hence, $t\notin V(T_{t'})$ and
        $t'\notin V(T_t)$, so $r$ separates $t$ from $t'$ in $T$, and so $V(A_{t'})\subseteq V(B_t)$. However,
        $V(A_{t'})\not\subseteq V(G_1)$, or else
        $V(A_{t'})\subseteq V(A_t\cup G_1)$, so, because $A_{t'}$ nests in
        $(G_1,G_2)$, $V(A_{t'})\subseteq V(G_2)$. But then
        $V(A_{t'})\subseteq V(B_t\cap G_2)$, contradiction.
      \end{subproof}
    \end{claim}

    Let $T_1$ and $T_2$ be disjoint copies of the tree $T$ and for each node
    $w\in V(T)$ let $w_1$ and $w_2$ denote the copies of $w$ in $T_1$ and $T_2$
    respectively. For each labelled leaf $t$ of $T$ we let $t_1=t$ if
    $t\in E(G_1)$ and we let $t_2=t$ if $t\in E(G_2)$. Now we construct a
    tree-decomposition $T'$ of $G$ by adding an edge $e_r=r_1r_2$ to
    $T_1\cup T_2$.

    \begin{claim}
      \label{thm:conn-set-decomposition-3}
      The adhesion of $T'$ is less than $\theta$.
      \begin{subproof}
        The edge $e_r$ has adhesion $\lambda(G_1,G_2)<\theta$. By symmetry it
        suffices to consider edges of $T_1$ in $T'$. Let $t\in V(T)-\{r\}$ and
        let $e'$ denote the copy of the edge $e_t$ in $T_1$. Note that the
        adhesion of $e'$ is at most $\lambda(G_1\cap A_t,G_2\cup B_t)$ which, by
        \autoref{thm:conn-set-decomposition-1}, is at most
        $\lambda(A_t,B_t)<\theta$.
      \end{subproof}
    \end{claim}

    \begin{claim}
      \label{thm:conn-set-decomposition-4}
      For each $t\in V(T)-\{r\}$ and each $i\in\{1,2\}$, we have
      $|Y_{t_i}(T')|\leq|Y_t(T)|$; moreover, if equality holds then
      $|Y_{t_{j}}(T')|\le \lambda(G_1,G_2)$ for some $j\in\{1,2\}$.
      \begin{subproof}
        By symmetry we may assume that $i=1$. Note that
        \[\begin{aligned}
            Y_{t_1}(T')-Y_t(T)&\subseteq V(A_t)\cap V(G_1) \cap V(G_2) - V(B_t),\\
            Y_t(T)-Y_{t_1}(T')&\supseteq V(A_t)\cap V(B_t)-V(G_1),
          \end{aligned} \]
        Then
        \[\begin{aligned}
            |Y_{t_1}&(T')| - |Y_t(T)| \\
            =&|Y_{t_1}(T')-Y_t(T)| - |Y_t(T) - Y_{t_1}(T')| \\
            \leq& |V(A_t)\cap V(G_1\cap V(G_2)-V(B_t)| - |V(A_t)\cap V(B_t)-V(G_1)|\\
            =& |V(A_t)\cap V(G_1)\cap V(G_2)-V(B_t)| \\
            &- |V(A_t)\cap V(B_t)| + |V(A_t)\cap V(B_t)\cap V(G_1)|\\
            =&|V(A_t)\cap V(G_1)\cap (V(B_t)\cup V(G_2))| - |V(A_t)\cap V(B_t)|\\
            =&\lambda(A_t\cap G_1,B_t\cup G_2) - \lambda(A_t,B_t)
          \end{aligned} \]
        so, by \autoref{thm:conn-set-decomposition-1},
        $|Y_{t_1}(T')|\leq|Y_t(T)|$. Moreover, if $|Y_{t_1}(T')|=|Y_t(T)|$,
        then, by \autoref{thm:conn-set-decomposition-1}, there exists
        $j\in\{1,2\}$ such that $V(A_t)\subseteq V(G_j)$, which implies that
        $Y_{t_{j+1}}(T')\subseteq V(G_1)\cap V(G_2)$, and hence
        $|Y_{t_{j+1}}|\leq\lambda(G_1,G_2)$, as required.
      \end{subproof}
    \end{claim}

    \begin{claim}
      \label{thm:conn-set-decomposition-5}
      For each $i\in\{1,2\}$, we have $|Y_{r_i}(T')|<|Y_{r}(T)|$.
      \begin{subproof}
        By symmetry we may assume that $i=1$. By our choice of $r$ and
        $(G_1,G_2)$ we have $|Y_{r}(T)\cap V(G_2)|>|V(G_1)\cap V(G_2)|$.
        Therefore,
        \[\begin{aligned}
            |Y_{r_1}(T')-Y_r(T)|\leq& |V(G_1)\cap V(G_2)-Y_r(T)| \\
            =& |V(G_1)\cap V(G_2)|-|Y_r(T)\cap V(G_1)\cap V(G_2)| \\
            <& |Y_r(T)\cap V(G_2)|-|Y_r(T)\cap V(G_1)\cap V(G_2)| \\
            =& |Y_r(T)-V(G_1)| \\
            \leq& |Y_r(T)-Y_{r_1}(T')|.
          \end{aligned} \]
        Therefore $|Y_{r_1}(T')|<|Y_{r}(T)|$.
      \end{subproof}
    \end{claim}

    By
    Claims~\ref{thm:conn-set-decomposition-1},~\ref{thm:conn-set-decomposition-4},
    and~\ref{thm:conn-set-decomposition-5}, the tree-decomposition $T'$ of $G$
    has adhesion less than $\theta$ and width at most $\omega$. By
    \autoref{thm:conn-set-decomposition-5}, we have
    $\max\{|Y_{r_1}(T')|, |Y_{r_2}(T')|\}< |Y_r(T)|$. Let
    $k\in\{0,\ldots,\omega+1\}$ be maximum such that there exists $s\in V(T)$
    with
    \[\max\{|Y_{s_1}(T')|,|Y_{s_2}(T')|\}<|Y_{s}|=k;\]
    thus $k \ge |Y_r| > \lambda(G_1,G_2)$. Consider any node $t\in V(T)$ with
    $|Y_t(T)|\ge k$ and $\max\{|Y_{t_1}(T')|,|Y_{t_2}(T')|\}=|Y_{t}|$. By
    \autoref{thm:conn-set-decomposition-4},
    $$\min\{|Y_{t_1}(T')|,|Y_{t_2}(T')|\}\le \lambda(G_1,G_2) < k.$$
    Therefore $n_k(T')<n_k(T)$ and $n_l(T') \le n_l(T)$ for each $l>k$. Thus
    $(n_{\omega+1}(T'),\ldots,n_0(T'))$ lexicographically precedes
    $(n_{\omega+1}(T),\ldots,n_0(T))$, contradicting our choice of $T$.
  \end{proof}
\end{theorem}

For $\theta\in \{1,2\}$, it is an easy exercise to show that the maximum size of
a $\theta$-connected set in a graph $G$ is equal to $\tw_{\theta}(G) +1$. Next
we will show, for $\theta\ge 3$, if a graph $G$ has a sufficiently large
$\theta$-connected set, then $G$ has large $\theta$-tree-width.

\begin{lemma}
  \label{thm:k-conn-set-to-tw}
  For integers $\theta$ and $n$ with $n\geq \theta\geq 3$, if $G$ is a graph
  with a $\theta$-connected set of size greater than
  $\binom{n}{\theta-1}(\theta-1)$, then $G$ has $\theta$-tree-width at least
  $n$.
  \begin{proof}
    Let $U$ be a $\theta$-connected set with $|U|>\binom{n}{\theta-1}(\theta-1)$
    and let $T$ be a tree-decomposition of $G$ with adhesion less than $\theta$
    and minimum width. We may choose $T$ so that all leaves are labelled and so
    that it has no vertices of degree two, and, subject to that, having as many
    vertices as possible. This choice of $T$ ensures that if $t$ is a vertex of
    degree at least $4$ in $G$ and $e$ and $e'$ are incident with $t$, then
    $Y_e\neq Y_{e'}$.

    Consider each edge $e=tt'$ of $T$ in turn. Let $G_{t}$ and $G_{t'}$ be the
    subgraphs of $G$ induced by the edges labelling the components of $T-e$ that
    contain the vertices $t$ and $t'$ respectively. Since $(G_t,G_{t'})$ is a
    $(\theta-1)$-separation of $G$ either $|V(G_{t})\cap U|<\theta$ or
    $|V(G_{t'})\cap U|<\theta$; in fact, because $|U|\geq 2\theta$, exactly one
    of $V(G_t)\cap U$ or $V(G_{t'})\cap U$ has size $\ge \theta$; orient the
    edge $e$ of $T$ toward the larger of the two sets.

    Since this oriented tree has no directed circuits, it has a vertex, say $r$,
    whose incident edges are all oriented toward it. Let $E_r$ denote the set of
    edges incident with $r$. Since
    $$|E_r|(\theta-1)\ge|U|> \binom{n}{\theta-1}(\theta-1),$$
    we have $|E_r|> \binom{n}{\theta-1}.$ In particular, $r$ has degree at least
    $4$. Note that $Y_{r}=\bigcup ( Y_e\, : \, e\in E_r)$ and, by our choice of
    $T$, the sets $(Y_e\, : \, e\in E_r)$ are pairwise distinct. Therefore
    $|E_r|\le \binom{|Y_{r}|}{\theta-1}$ and hence $|Y_{r}|\ge n+1$. So
    $\tw_\theta(G)\ge n$.
  \end{proof}
\end{lemma}

Combining Theorems~\ref{thm:conn-set-decomposition} and~
\ref{thm:k-conn-set-to-tw}, we get
\autoref{thm:k-tree-width-connected-set-duality}.

The following construction shows that the bound $\binom{n}{\theta-1}(\theta-1)$
in \autoref{thm:k-conn-set-to-tw} is tight:
Begin with an empty graph on a vertex set $A$ with $|A|=n$ and let $\cA$ denote
the set of all $A'\subseteq A$ with $|A'|=\theta-1$. For each set $A'\in \cA$,
add a set $B_{A'}$ of $\theta-1$ new vertices and add an edge between each
$a\in A'$ and $b\in B_{A'}$. Let $B$ be the union of the sets $B_{A'}$, over all
sets $A'\in\cA$. Then $\tw_\theta(G)\leq n-1$ and $B$ is a $\theta$-connected set of
size $\binom{n}{\theta-1}(\theta-1)$.

\section{Necklaces}

\begin{figure}[]
  \centering
  \includesvg[clean,width=0.5\columnwidth]{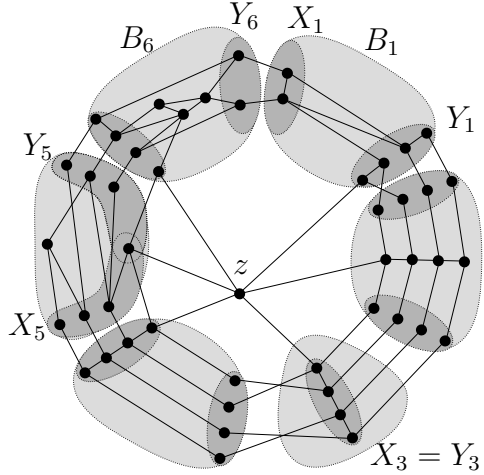}
  \caption{An example of a $(4,2,1,6)$-necklace. }
  \label{fig:necklace-defn}
\end{figure}

Directly constructing generalized-wheel-minors is unwieldy, so we start by
constructing a closely related structure, called a \emph{necklace}; see
\autoref{fig:necklace-defn}.

Let $t$, $s$, $\ell$ and, $n$ be non-negative integers with $t\geq s$ and
$n\geq 2$. We denote the additive group on the integers modulo $n$ by $\Z_n$.
Let $G$ be a graph and let $N=(\cB,\cM,Y)$ where $\cB$ is a sequence
$(B_i\, :\, i\in \Z_n)$ of subgraphs of $G$, $\cM=(M_i\, :\, i\in\Z_n)$ is a
sequence of edge sets in $G$, and $Z\subseteq V(G)$. The {\em left attachment
  sequence} of $N$ is the sequence $\cX=(X_i\, :\, i \in \Z_n)$ where $X_i$
denotes the set of all vertices in $B_i$ that are incident with an edge in
$M_{i-1}$, for each $i\in \Z_n$. The {\em right attachment sequence} of $N$ is
the sequence $\cY=(Y_i\, :\, i \in \Z_n)$ where $Y_i$ denotes the set of all
vertices in $B_i$ that are incident with an edge in $M_{i}$, for each
$i\in \Z_n$. We call $N$ a {\em $(t,s,\ell,n)$-necklace} in $G$ if it satisfies
the following properties:
\begin{enumerate}[(N$_1$)]
\item the subgraphs $(B_i\, :\, i\in \Z_n)$ are pairwise disjoint and are each
  disjoint from $Z$,
\item for each $i\in \Z_n$, $B_i$ is connected and non-empty,
\item for each $i\in \Z_n$, the set $M_i$ is a matching where each edge has one
  end in $B_i$ and one end in $B_{i+1}$,
\item each of $M_1,\ldots,M_{n-1}$ has size $t$ and $M_n$ has size $s$,
\item for each $i\in\{2,\ldots,n-1\}$, there are $t$ vertex-disjoint
  $(X_i,Y_i)$-paths in $B_i$,
\item for both $i\in\{1,n\}$, there are $s$ vertex-disjoint $(X_i,Y_i)$-paths in
  $B_i$, and
\item $|Z|=\ell$ and every vertex in $Z$ has a neighbour in each of the sets
  $V(B_1),\ldots,V(B_n)$.
\end{enumerate}
We call the graphs in $\cB$ the \emph{beads} and the vertices in $Z$ the
\emph{hubs} of the necklace. By a {\em $(\theta;n)$-necklace} we mean a
$(t,s,\ell,n)$-necklace where $t+s+\ell=\theta$.

We say that a $(t,s,\ell,n)$-necklace is {\em supported by} a set
$U\subseteq V(G)$ if each bead contains at least one vertex in $U$. The
following theorem is the main step toward proving our main results,
\autoref{thm:nearly-balanced-necklace} and \autoref{thm:main-local}.

\begin{theorem}
  \label{thm:set-to-necklace-preview}
  For all non-negative integers $\theta$ and $n$ with $n\ge 2$ there is an
  integer $m$ such that, if $G$ is a graph containing a $\theta$-connected set
  $U$ of size at least $m$, then $G$ contains a $(\theta;n)$-necklace supported
  by $U$.
\end{theorem}

The proof of \autoref{thm:set-to-necklace-preview} essentially goes by induction
on $\theta$. In \autoref{sec:init} we find a $(0,0,\theta,n)$-necklace or a
$(1,0,0,m)$-necklace with $m\gg n$. Then, in \autoref{sec:jump} and
\autoref{sec:jump-free}, we show how to turn the $(1,0,0,m)$-necklace into a
$(\theta;n)$-necklace.

We conclude this section by defining notation that will be used in subsequent
sections.

For brevity, when we say that $N$ is a $(t,s,\ell,n)$-necklace we implicitly
assume that $N=(\cB,\cM,Y)$ where $\cB=(B_i\, :\, i\in \Z_n)$ and
$\cM=(M_i\, :\, i\in \Z_n)$, and moreover, we assume that the left and right
attachment sequences are $\cX=(X_i\, :\, i \in \Z_n)$ and
$\cY=(Y_i\, :\, i \in \Z_n)$. Similarly, when $N'$ is a $(t',s',\ell',n')$-necklace we
implicitly assume that $N'=(\cB',\cM',Z')$ where $\cB'=(B'_i\, :\, i\in \Z_{n'})$
and $\cM'=(M'_i\, :\, i\in \Z_{n'})$, and that the left and right attachment
sequences are $\cX'=(X'_i\, :\, i \in \Z_{n'})$ and $\cY'=(Y'_i\, :\, i \in \Z_{n'})$.
We make similar assumptions for necklaces $N''$, $N^1$, $N^2$, et cetera.

Let $N$ be a $(t,s,\ell,n)$-necklace in $G$. We denote
$V(B_1\cup\cdots\cup B_n)\cup Z$ by $V(N)$. We let $G[N]$ denote the  subgraph of
$G$ whose vertex set is $V(N)$ and whose edge set is the union of
$E(B_1\cup\cdots\cup B_n)$, $M_1\cup\cdots\cup M_n$, and the set of all edges
with one end in $Z$ and the other end in $V(B_1\cup\cdots\cup B_n)$.

The \emph{reversal} of $N$ is the $(t,s,\ell,n)$-necklace $N'=(\cB',\cM',Z)$
where, for each $i\in\Z_n$, we have $B'_i=B_{n-i+1}$ and $M'_i =M_{n-i}$.

We say that $N$ {\em supports} a necklace $N'$ if each bead of $N'$ contains a
bead of $N$ as a subgraph. Note that this is a transitive relation on necklaces
and that, if $N'$ is supported by $N$ and $N$ is supported by a vertex set $U$,
then $N'$ is supported by $U$ as well.

Suppose $N$ is a $(t,s,\ell,m)$-necklace in a graph $G$. Let $(a_0,\ldots,a_n)$
be a sequence of numbers with $0=a_0<\cdots<a_n=m$. For each $i\in \Z_n$, let
$B'_i$ denote the subgraph of $G$ obtained taking the union of the beads
$(B_{a_{i-1}+1},B_{a_{i-1}+2},\ldots,B_{a_{i}})$ together with the matchings
$(M_{a_{i-1}+1},M_{a_{i-1}+2},\ldots,M_{a_{i}-1})$. For each
$i\in\{1,\ldots,n\}$, let $M'_i = M_{a_i}$. Finally let $\cB'$ and $\cM'$ denote
the sequences $(B'_i\,:\, i\in \Z_n)$ and $(M'_i\,:\, i\in \Z_n)$, respectively. We
call $N'=(\cB',\cM',Z)$ the \emph{contraction of $N$ to $(a_1,\ldots,a_{n-1})$}
and we denote it by $N\circ (a_1,\ldots,a_{n-1})$. Note that, if $t\neq 0$, then
$N'$ is a $(t,s,\ell,n)$-necklace in $G$ supported by $N$.

\section{Initial necklace}
\label{sec:init}

In this section we prove the base case for the inductive proof of
\autoref{thm:set-to-necklace-preview}. \autoref{lem:init-necklace} shows that a
sufficiently large $\theta$-connected set gives rise to either a
$(0,0,\theta,n)$-necklace or a $(1,0,0,m)$-necklace with $m\gg n$.

A connected graph $G$ with $|V(G)|\ge d^n$ must contain either an $n$-vertex
path or a vertex of degree $d$; the following result is just a refinement of
that observation. The bound of $d^n$ is larger than necessary, but this simplifies
the proof slightly.

\begin{lemma}
  \label{lem:long-path-or-high-degree}
  For positive integers $d$, $n$ and $m$ with $d\geq 2$, $n\geq 2$ and
  $m\geq d^n$, if $N$ is a $(0,0,0,m)$-necklace in a connected graph $G$, then
  either
  \begin{enumerate}[(i)]
  \item there is a $(0,0,1,d)$-necklace in $G$ that is supported by $N$, or
  \item there is a $(1,0,0,n)$-necklace in $G$ that is supported by $N$.
  \end{enumerate}
  \begin{proof}
    Consider a counterexample $(d,n,m,N,G)$ with $|E(G)|$ minimum. By growing
    each of the beads we may assume that $V(G) = V(N)$. By the minimality of
    $G$,
    \begin{enumerate}[(a)]
    \item $B_1,\ldots,B_m$ are trees,
    \item $G$ is a tree, and
    \item if $v$ is a vertex of $G$ with degree at most two, then $v$ induces a
      singleton bead in $\cB$ (since otherwise we would obtain a smaller
      counterexample by contracting an edge of $B_1\cup\cdots\cup B_m$ that is
      incident with $v$).
    \end{enumerate}

    First suppose that $G$ has a vertex $z$ with degree at least $d$. Let
    $(B'_1,\ldots,B'_{d})$ be a sequence of distinct components of $G-v$. Each
    of $B'_1,\ldots,B'_{d}$ has a vertex that has degree at most one in $G$, so,
    by $(c)$, each of $B'_1,\ldots,B'_{d}$ contains a bead of $N$. Thus
    $(B'_1,\ldots,B'_{d})$ is the bead sequence of a $(0,0,1,d)$-necklace
    supported by $N$. Thus $(i)$ holds.

    Suppose now that each vertex in $G$ has degree less than $d$. Since
    $|V(G)|>d^n$, $G$ contains an $n$-vertex path $P$.
    Each component of $G\delete E(P)$
    contains a vertex whose degree in $G$ is at most two, so each component of
    $G\delete E(P)$ contains a bead of $N$ by $(c)$. These components together
    with $E(P)$ form a $(1,0,0,n)$-necklace in $G$, so $(ii)$ holds.
  \end{proof}
\end{lemma}

A simple inductive argument now gives a $(0,0,\theta,n)$- or
$(1,0,0,p)$-necklace from a sufficiently large $\theta$-connected set.

\begin{lemma}
  \label{lem:init-necklace}
  There exists a function $\fcur:\N^3\rightarrow\N$ such that, for
  $\theta,n,p\in\N$, with $n\geq 2$ and $p\geq 2$, if $U$ is a
  $\theta$-connected set in a graph $G$ with $|U|\ge \fcur(\theta,n,p)$, then
  either
  \begin{enumerate}[(i)]
  \item $G$ contains a $(0,0,\theta,n)$-necklace supported by $U$, or
  \item $G$ contains a $(1,0,0,p)$-necklace whose beads each contain at least
    $\theta$ vertices in $U$.
  \end{enumerate}
  \begin{proof}
    Fix $\theta,n,p\in\N$. Let $m_\theta=n$ and, for
    $k\in\{0,\ldots,\theta-1\}$, let $m_k=m_{k+1}^{p\theta}+k$. Let
    $\fcur(\theta,n,p)=m_0$. Let $U$ be a $\theta$-connected set in a graph $G$
    with $|U|\geq \fcur(\theta,n,p)$. Note that the singletons of elements in
    $U$ form a $(0,0,0,m_0)$-necklace; let $k\in\{0,\ldots,\theta\}$ be maximal
    such $G$ contains a $(0,0,k,m_k)$-necklace $N$ supported by $U$. If
    $k=\theta$, then $(i)$ holds.

    Now suppose $k<\theta$. Let $Z$ be the hub set of $N$, so $|Z|=k$. Because
    $U$ is $\theta$-connected, $G-Z$ contains a component $G'$ such that
    $|U-V(G')|\leq k$. Therefore, $|U\cap V(G')|\geq m_{k+1}^{ p\theta}$. By
    \autoref{lem:long-path-or-high-degree}, $G'$ contains either a
    $(0,0,1,m_{k+1})$- or $(1,0,0,p \theta)$-necklace $N'$ supported by $N$.
    Each $z\in Z$ is a hub of $N'$, so if $N'$ is a $(0,0,1,m_{k+1})$-necklace,
    then $G$ contains a $(0,0,k+1,m_{k+1})$-necklace supported by $N$, and hence
    supported by $U$, contradicting our choice of $k$. Therefore, $N'$ is a
    $(1,0,0,p \theta)$-necklace. Thus,
    $N'\circ (\theta,2\theta,\ldots,(p-1)\theta)$ is a $(1,0,0,p)$-necklace in
    $G$ with each bead containing at least $\theta$ vertices in $U$, so $(ii)$
    holds.
  \end{proof}
\end{lemma}

\section{Long jumps}
\label{sec:jump}

In the next two sections we describe how to turn a $(1,0,0,m)$-necklace, with
$m\gg n$, into a $(\theta;n)$-necklace.

Let $N$ be a $(t,s,\ell,n)$-necklace in a graph $G$. For $i,j\in\Z_n$, an
\emph{$(i,j)$-jump} of $N$ is a $(V(B_i),V(B_j))$-path $P$ which intersects
$V(N)$ only at its ends.
A {\em long jump} is an
$(i,j)$-jump where $j-i\not\in\{-1,0,1\}$; note that a $(1,n)$-jump is not a
long jump. If $G$ does not contain a long jump of $N$, then we say that $N$ is
{\em long-jump-free} in $G$.

In this section we turn a $(t,s,\ell,m)$-necklace, with $t\geq 1$ and $m\gg n$,
into either a long-jump-free $(t,s,\ell,n)$-necklace or a
$(t',s',\ell',n)$-necklace with $t'\geq 1$ and $t'+s'+\ell' > t+s+\ell$.

\begin{lemma}
  \label{lem:jump-structure}
  There exists a function $\fcur:\N\rightarrow\N$ such that the following holds.
  Let $t,s,\ell,m,n\in\N$ with $t\geq 1$, $n\geq 2$ and $m\geq\fcur(n)$. If $G$
  is a graph containing a $(t,s,\ell,m)$-necklace $N$, then $G$ contains a
  $(t,s,\ell,n)$-necklace $N'$ supported by $N$ such that either
  \begin{enumerate}[(i)]
  \item $N'$ is long-jump-free in $G$, or
  \item for each $i\in\Z_{n}$, there is a $(1,i)$-jump for $N'$ in $G$.
  \end{enumerate}
  \begin{proof}
    Fix $n\in\N$. For $k\in\{0,\ldots,n^2+1\}$, let
    $m_k=n^{n^2-k+1}$. Let $\fcur(n+2)=m_0+2$. Let $t,s,\ell,m\in\N$ with
    $t\geq s$ and $m\geq \fcur(n+2)$, and let $G$ be a graph containing a
    $(t,s,\ell,m)$-necklace $N$. It suffices to show that $G$ contains a
    $(t,s,\ell,n+2)$-necklace $N'$ supported by $N$ such that either $(i)$ $N'$
    is long-jump-free or $(ii)$ for each $i\in\Z_{n+2}$, there is a $(1,i)$-jump
    for $N'$ in $G$.
  
    Let $S$ denote the set of all pairs $(a,b)$ such that $1\le a\leq b\le m$
    and there exists an $(a,b)$-jump for $N$ in $G$. For integers $a$ and $b$
    with $a\leq b$ we let $[a,b]$ denote the interval $\{a,a+1,\ldots,b\}$. A
    subset $S'\subseteq S$ is {\em $t$-intersecting} if
    \[\left|[1,m]\cap\bigcap_{(a,b)\in S'}[a,b]\right|\geq t.\]

    Let $S'\subseteq S$ be maximal such that $S'$ is $(m_k+2)$-intersecting,
    where $k=|S'|$; note that $S'$ exists because $|[1,m]|\geq m_0+2$, so
    $\emptyset$ is $(m_0+2)$-intersecting. Let
    $\{(a_1,b_1),\ldots,(a_k,b_k)\}=S'$. Let $a_{\max}=\max\{1,a_1,\ldots,a_k\}$,
    $b_{\min}=\min\{m,b_1,\ldots,b_k\}$ and let $I=[a_{\max},b_{\min}]$, and note
    that
    \[a_{\max}+m_k< a_{\max}+|I|-1=b_{\min}.\]

    \begin{case-one}
      \quad $k\leq n^2$.
    \end{case-one}

    Then $a_{\max}+n\cdot m_{k+1}<b_{\min}$. Let
    $$N'= N\circ(a_{\max}, a_{\max}+m_{k+1},\ldots, a_{\max}+n\cdot m_{k+1}).$$
    Thus $N'$ is a $(t,s,\ell,n+2)$-necklace supported by $N$. By our choice of
    $S$, for each $(a,b)\in S-S'$, $|I\cap [a,b]|\leq m_{k+1}+1$, so
    $V(B_a)\subseteq V(B'_{a'})$ and $V(B_b)\subseteq V(B'_{b'})$ for some $a'$
    and $b'$ with $a'\leq b' \leq a'+1$; for each $(a,b)\in S'$,
    $V(B_a)\subseteq V(B'_1)$ and $V(B_b)\subseteq V(B'_{n+2})$. Therefore $N'$
    is long-jump-free in $G$ and, hence, $(i)$ holds.

    \begin{case-two}
      \quad $k=n^2+1$.
    \end{case-two}

    Then $\max\{|\{a_1,\ldots,a_k\}|,|\{b_1,\ldots,b_k\}|\}> n$, so by
    possibly reversing $N$, we may assume that $|\{b_1,\ldots,b_k\}|\ge n+1$. By
    possibly relabeling the elements of $S'$, we may assume that
    $b_1<b_2<\cdots< b_{n+1}$. Also note that
    $a_{\max}<a_{\max}+m_k<b_{\min}\leq b_1$. Let
    $$N'= N\circ (a_{\max},b_1,b_2,\ldots,b_{n}).$$
    Then $N'$ is a $(t,s,\ell,n+2)$-necklace and, for each
    $i\in\{2,\ldots,n+2\}$, the $(a_{i-1},b_{i-1})$-jump in $N$ is a
    $(1,i)$-jump in $N'$, so $(ii)$ holds.
  \end{proof}
\end{lemma}

We will now treat outcome $(ii)$ of \autoref{lem:jump-structure}. For
convenience we will handle the trivial case where $s=0$ in the following lemma
whose proof we omit.
\begin{lemma}
  \label{lem:s-zero}
  Let $N$ be a $(t,0,\ell,n)$-necklace in a graph $G$ with $t\geq 1$. If there
  is a $(1,n)$-jump of $N$ in $G$, then $G$ has a $(t,1,\ell,n)$-necklace
  supported by $N$.
\end{lemma}

\begin{lemma}
  \label{lem:grow-jumps-from-weak-bead}
  There exists a function $\fcur:\N^2\rightarrow\N$ such that the following
  holds. Let $t,s,m,n\in\N$ with $t\ge s\ge 1$, $n\ge 2$, and
  $m\geq \fcur(s,n)$. If $N$ is a $(t,s,\ell,m)$-necklace in a graph $G$ and
  there is a $(1,i)$-jump of $N$ for each $i\in\{1,\ldots,m-1\}$, then there
  exists a $(t+1,s-1,\ell,n)$- or $(t,s-1,\ell+1,n)$-necklace supported by $N$
  in $G$.
  \begin{proof}
    Fix $s,n\in\N$ with $n\geq 2$ and $s\geq 1$. Let
    \[\fcur(s,n) = \left( n^{(n-1)^2+1}-1 \right)s+4.\]
    Let $t,m\in\N$ with $t\geq s$ and $m\geq \fcur(s,n)$. Let $N$ be a
    $(t,s,\ell,m)$-necklace in a graph $G$ such that, for each
    $i\in\{1,\ldots,m-1\}$, $G$ contains a $(1,i)$-jump for $N$.
  
    By possibly growing the bead $B_1$ we may assume that each $(1,i)$-jump
    consists of a single edge. For each $i\in \{3,4,\ldots,m-1\}$ let $e_i=x_iy_i$
    be an edge with $x_i\in V(B_1)$ and $y_i\in V(B_i)$. There is a collection
    of vertex-disjoint trees $(T_1,\ldots, T_s)$ in $B_1$ such that
    $V(B_1)=V(T_1\cup\cdots \cup T_s)$ and, for each $i\in\{1,\ldots,s\}$, the
    tree $T_i$ contains a vertex in each of $X_1$ and $Y_1$. Then
    $|\{3,\ldots,m-1\}|>\left( n^{(n-1)^2+1}-1 \right)s$, so there is an integer
    $j\in\{1,\ldots,s\}$ and a set $I\subseteq \{3,\ldots,m-1\}$ with
    $|I|=n^{(n-1)^2+1}$ such that $\{x_i\, : \, i\in I\}\subseteq V(T_j)$.
    Without loss of generality, $j=s$. Let $T$ be the tree obtained from $T_s$
    by adding the vertices $\{y_i\, : \, i\in I\}$ and the edges
    $\{e_i\, : \, i\in I\}$.

    \begin{claim}
      \label{lem:grow-jumps-from-weak-bead:claim-N2}
      There exists a $(0,0,1,n)$- or $(1,0,0,n)$-necklace $N^2$ in $T$ and an
      increasing sequence $(a_1,a_2,\ldots,a_n)$ of integers in $I$ such that,
      for $i\in\{1,\ldots,n\}$, $y_{a_i}\in B^2_i$.
      \begin{subproof}
        By \autoref{lem:long-path-or-high-degree}, there is a $(0,0,1,n)$- or
        $(1,0,0,(n-1)^2+1)$-necklace $N^1$ in $T$ that is supported by
        $\{y_i\, : \, i\in I\}$.

        Suppose that $N^1$ is a $(0,0,1,n)$-necklace.
        Then any permutation of the beads of $N^1$ yields another
        $(0,0,1,n)$-necklace. In particular, the beads can be ordered such that
        $(a_1,\ldots,a_n)$ exists.

        Now suppose that $N^1$ is a $(1,0,0,(n-1)^2+1)$-necklace. For
        $i\in\{1,\ldots,(n-1)^2+1\}$, choose $b_i\in I$ such that
        $y_{b_i}\in V(B^1_i)$. By the Erd{\H o}s-Szekeres Theorem,
        $(b_1,\ldots,b_{(n-1)^2+1})$ contains a monotonic subsequence
        $(a_1,\ldots,a_n)$. By possibly reversing $N$, we may assume that
        $a_1<\cdots<a_n$. By contracting $N^1$, we obtain a $(1,0,0,n)$-necklace
        $N^2$ such that, for $i\in\{1,\ldots,n\}$, $y_{a_i}\in V(B^2_i)$.
      \end{subproof}
    \end{claim}

    Let $N^3 = N\circ (a_1,\ldots,a_{n-1})$. Note that, because $a_1\ge 3$,
    $B^3_1$ contains both $B_1$ and $B_2$ and, hence, $B^3_1-V(T_s)$ is
    connected. For $j\in\{1,\ldots,s-1\}$, $T_j$ contains an $(X_1,Y_1)$-path, so
    $B^3_1-V(T_s)$ contains $s-1$ vertex-disjoint $(X_1,Y^3_1)$ paths. Also note
    that $M^3_n$ contains exactly one edge, $e_s$, incident with a vertex in
    $V(T_s)$. So we can construct a $(t,s-1,\ell,n)$-necklace $N^4$ by deleting
    $V(T_s)$ from $B^3_1$, and deleting $e_s$ from $M^3_n$. Note that $N^4$ is
    supported by $N$ because $B_2\subseteq B^4_1$.
    See \autoref{fig:many-jumps}.
    \begin{figure}
      \centering
      \includesvg[clean,width=0.5\columnwidth]{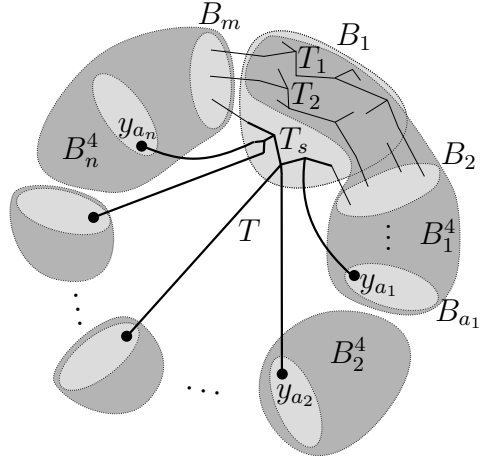}
      \caption{Notation used in the proof of
        \autoref{lem:grow-jumps-from-weak-bead}. The beads of $N$ are shown in
        light grey, and the beads of $N^4$ are shown in dark grey.}
      \label{fig:many-jumps}
    \end{figure}

    Note that, for $i\in\Z_n$,
    $V(B^2_i)\cap V(N^4)=\{y_{a_i}\}\subseteq V(B^4_i)$,
    $Z^2\cap V(N^4)=\emptyset$, and $M^2_i\cap E(N^4)=\emptyset$. We construct a
    $(t,s-1,\ell+1,n)$- or $(t+1,s-1,\ell,n)$-necklace $N^5$ by merging $N^2$
    with $N^4$ as follows:
    \begin{enumerate}[(i)]
    \item let $Z^5=Z^2\cup Z^4$,
    \item for each $i\in\Z_n$, let $B^5_i=B^2_i\cup B^4_i$, and
    \item for each $i\in \Z_n$, let $M^5_i=M^2_i\cup M^4_i$.
    \end{enumerate}
    Finally, note that $N^5$ is supported by $N$ via $N^4$, as required.
  \end{proof}
\end{lemma}

Finally combining Lemmas~\ref{lem:jump-structure}, \ref{lem:s-zero},
and~\ref{lem:grow-jumps-from-weak-bead} with a routine inductive proof on $s$
(which is omitted) we obtain the following result.

\begin{lemma}
  \label{lem:long-jump-free}
  There exists a function $\fcur:\N^2\rightarrow\N$ such that the following
  holds. Let $t,s,m,n,\ell \in\N$ with $t\ge 1$, $n\geq 2$ and
  $m\geq\fcur(s,n)$. If $N$ is a $(t,s,\ell,m)$-necklace in a graph $G$, then
  there exists a $(t',s',\ell',n)$-necklace $N'$ in $G$ supported by $N$ such
  that $t'\geq 1$, and either
  \begin{enumerate}[(i)]
  \item $t'+s'+\ell'=t+s+\ell$ and $N'$ is long-jump-free in $G$, or
  \item $t'+s'+\ell'>t+s+\ell$.
  \end{enumerate}

\end{lemma}

\section{Long-jump-free necklaces}
\label{sec:jump-free}
Henceforth we may assume that we have a long-jump-free $(t,s,\ell,n)$-necklace,
supported by a $\theta$-connected set $U$, in a graph $G$, where
$t+s+\ell <\theta$. In this section we show that $G$ contains
$(t+1,s,\ell,n')$-necklace or a $(t,s+1,\ell,n')$-necklace
(\autoref{lem:grow-theta}), and then complete the proof of
\autoref{thm:set-to-necklace-preview}, as \autoref{thm:necklace-main}.

First, we require the following easy connectivity result.
\begin{lemma}
  \label{lem:connectivity}
  Let $k\in\N$. Let $\cP$ be a collection of $k$ vertex-disjoint $(X,X')$-paths
  in a graph $G$, and let $\cQ$ be a collection of $k$ vertex-disjoint
  $(Y,Y')$-paths in $G$, where $X,X',Y,Y'\subseteq V(G)$. Let $\cB$ be a
  collection of $k$ vertex-disjoint connected subgraphs of a graph $G$, each
  having non-empty intersection with each path in $\cP\cup\cQ$. Then $G$
  contains $k$ vertex-disjoint $(X,Y)$-paths.

  \begin{proof}
    Consider any set $A\subseteq V(G)$ with $|A|<k$. Since
    $|\cB|=|\cP|=|\cQ|= k>|A|$, there exist $B\in\cB$, $P\in \cP$, and
    $Q\in \cQ$ such that $A\cap V(P\cup B\cup Q)=\emptyset$. By the definition
    of $X$ and $Y$, we have $X\cap V(P)\neq \emptyset$ and
    $Y\cap V(Q)\neq \emptyset$. Moreover, since $B$ has a non-empty intersection
    with both $P$ and $Q$, the graph $B\cup P\cup Q$ is connected. So there
    exists an $(X,Y)$-path in $G-A$. Then, by Menger's Theorem, there are $k$
    vertex-disjoint $(X,Y)$-paths in $G$.
  \end{proof}
\end{lemma}

The following lemma shows how to exploit excess of paths through a segment of a
necklace to increase the structure. This is the key step in \autoref{lem:grow-theta}.
The proof of this result incorporates ideas 
from the short proof of the Grid Theorem due
to Diestel, Gorbunov, Jensen, and Thomassen~\cite{Diestel1999}.

\begin{lemma}
  \label{lem:grow-linear}
  There exists a function $\fcur:\N^3\rightarrow \N$ such that the following
  holds. Let $t,s,m,n\in\N$ with $t\geq s$, $n\geq 2$ and $m\geq \fcur(t,s,n)$. Let $G$ be a
  graph with a $(0,0,0,m)$-necklace $N$ that is long-jump-free and has no
  $(1,m)$-jump. Let $X\subseteq V(B_1)$ and $Y\subseteq V(B_m)$. Suppose $G$
  contains $s$ vertex-disjoint $(X,V(B_s))$-paths, $s$ vertex-disjoint
  $(V(B_{m-s+1}),Y)$-paths, and $t$ vertex-disjoint $(V(B_1),V(B_m))$-paths.
  Then $G$ contains a $(t,0,0,n)$-necklace $N'$ supported by $N$ such that
  $B'_1$ contains $s$ vertex-disjoint $(X,X'_1)$-paths and $B'_n$ contains $s$
  vertex-disjoint $(X'_n,Y)$-paths.
  \begin{proof}
    Fix $t,s,n\in\N$ with $t\geq s$ and $n\geq 2$. Let
    \[\fcur(t,s,n)=\max\{n\,,\,n(t-1)+2s-1\}.\]
    Let $N$ be a $(0,0,0,m)$-necklace in a graph $G$, with $m\geq\fcur(t,s,n)$,
    let $X\subseteq V(B_1)$, and let $Y\subseteq V(B_n)$. Suppose $G$ contains
    $s$ vertex-disjoint $(X,V(B_s))$-paths, $s$ vertex-disjoint
    $(V(B_{m-s+1}),Y)$-paths, and $t$ vertex-disjoint $(V(B_1),V(B_m))$-paths.
    If $t\leq 1$, then it easy to see that $G$ contains a $(t,0,0,n)$-necklace
    with the required properties. We now assume $t\geq 2$.

    \begin{claim}
      There exist sets $X'$ and $Y'$ such that $X'\subseteq V(B_1\cup B_{s+1})$,
      $Y'\subseteq V(B_{m-s}\cup B_m)$, $|X'\cap X|=|Y'\cap Y|=s$,
      $|X'|=|Y'|=t$, and $G$ contains $t$ vertex-disjoint
      $(X',Y')$-paths.
      \begin{subproof}
        Let $(P_1,\ldots,P_t)$ be a collection of vertex-disjoint
        $(V(B_1),V(B_m))$-paths, and let $\cP=(P_1,\ldots,P_{s})$. Let $\cQ$ be
        a collection of $s$ vertex-disjoint $(X,V(B_s))$-paths in $G$; because
        $N$ is long-jump-free and has no $(1,m)$-jump, we may assume each
        $Q\in\cQ$ is disjoint from $V(B_{s+1}\cup\cdots\cup B_m)$. Similarly,
        there exists a collection $\cQ'$ of $s$ vertex-disjoint
        $(V(B_{m-s+1}),Y)$-paths in $G$, each of which is disjoint from
        $V(B_{1}\cup\cdots\cup B_{m-s})$. Let $H$ be the union of the paths in
        $\cP\cup\cQ\cup\cQ'$ together with the beads $B_1,\ldots,B_s$ and
        $B_{m-s+1},\ldots,B_m$.

        Each path in $\cP\cup\cQ$ intersects each of $B_1,\ldots,B_s$, so by
        \autoref{lem:connectivity}, $H$ contains a collection $\cP'$ of $s$
        vertex-disjoint $(X,V(B_n))$-paths. Similarly, each path in
        $\cP'\cup\cQ'$ intersects each of $B_{m-s+1},\ldots,B_m$, so, by
        \autoref{lem:connectivity}, $H$ contains a collection
        $(P''_1,\ldots,P''_{s})$ of vertex-disjoint $(X,Y)$-paths.

        For $k\in\{s+1,\ldots,t\}$, $P_k-V(H)$ contains a
        $(V(B_{s+1}),V(B_{m-s}))$-path $P''_k$, because $s+1\leq m-s$ and $N$ is
        long-jump-free and has no $(1,m)$-jump. Let $x_k$ be the end of $P''_k$
        in $V(B_{s+1})$ and let $y_k$ be the end of $P''_k$ in $V(B_{m-s})$. Let
        $X'=X\cup\{x_{s+1},\ldots,x_t\}$ and let
        $Y'=Y\cup\{y_{s+1},\ldots,y_t\}$. Then, $(P''_1,\ldots,P''_t)$ is a
        collection of vertex-disjoint $(X',Y')$-paths, as desired.
      \end{subproof}
    \end{claim}

    Let $G'$ be the minimal subgraph of $G$ containing $G[N]$, and
    containing $t$ vertex-disjoint $(X',Y')$-paths. Let $(P_1,\ldots,P_{t})$
    be a collection of $t$ vertex-disjoint $(X',Y')$-paths. For
    $i\in\{s+2,\ldots,m-s\}$ and $k\in\{1,\ldots,t\}$, let $L_{i,k}$ denote
    the longest subpath of $P_k$ with an end in $X'$ and disjoint from $B_{i}$;
    let $R_{i,k}=P_{k}-V(L_{i,k})$. Note that both $L_{i,k}$ and $R_{i,k}$ are
    non-empty. Let $e_{i,k}$ be the edge in $P_k$ between $L_{i,k}$ and
    $R_{i,k}$. Let $L_{i}=\bigcup_{k=1}^{t}L_{i,k}$, and let
    $R_i=\bigcup_{k=1}^{t}R_{i,k}$.

    \begin{claim}
      \label{lem:grow-linear:backtrack}
      For each $i\in\{s+2,\ldots,m-s\}$ and each $j\in\{1,\ldots,m\}$ with
      $j\leq i-t+1$, $V(R_i\cap B_j)=\emptyset$.
      \begin{subproof}
    \begin{figure}
      \centering
      \includesvg[clean,width=0.5\columnwidth]{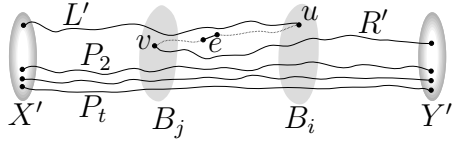}
      \caption{Notation used in the proof of
        \autoref{lem:grow-linear}.}
      \label{fig:jump-free-backtrack}
    \end{figure}
        Suppose the contrary, so there exists
        $v\in V(R_{i,k}\cap B_j)\neq\emptyset$, for some $k$; without loss of
        generality, $k=1$. Let $u$ be the end of $R_{i,1}$ in $B_i$. Note that
        $i\neq j$, so the subpath of $P_{1}$ between $u$ and $v$ contains some
        edge, $e$, that is not in $E(N)$. Let
        $L'$ denote the subpath of $P_{1}$ between $X'$ and $u$, and let $R'$
        denote the subpath of $P_{1}$ between $v$ and $Y'$. See \autoref{fig:jump-free-backtrack}.
        Then $(L',P_2,\ldots,P_t)$ and $(R',P_2,\ldots,P_t)$ are each
        collections of vertex-disjoint paths in $G'\delete e$ with ends in $X'$ and $Y'$, respectively.
        Moreover, each of the paths $P_2,\ldots,P_t$, $L'$, $R'$
        intersect each of $B_j,\ldots,B_i$, and $|\{B_j,\ldots,B_i\}|\geq t$.
        By \autoref{lem:connectivity}, 
        $G'\delete e$ contains $t$ vertex-disjoint $(X',Y')$-paths,
        contradicting the choice of $G'$.
      \end{subproof}
    \end{claim}

    For $i\in\{1,\ldots,n\}$, let $\alpha(i)=(i-1)(t-1)+s+1$. Note that
    $s+1=\alpha(1)<\cdots<\alpha(n)\leq m-s$.
    Let $B'_1=B_{\alpha(1)}\cup L_{\alpha(2)}$, let $B'_n=B_{\alpha(n)}\cup R_{\alpha(n)}$, and for
    $i\in\{2,\ldots,n-1\}$, let $B'_i=B_{\alpha(i)}\cup (L_{\alpha(i+1)}\cap R_{\alpha(i)})$.
    For $i\in\{1,\ldots,n-1\}$, let
    $M'_i=\{e_{\alpha(i+1),k}\,:\,k\in\{1,\ldots,t+1\}\}$. For
    $i\in\{1,\ldots,n-1\}$, let $Y'_i=V(B'_i\cap M'_i)$ and
    $X'_{i+1}=V(B'_{i+1}\cap M'_i)$.

    Note that $L_{\alpha(i+1)}$ is disjoint
    from $B_{\alpha(j)}$, for each $j>i$ and, by \autoref{lem:grow-linear:backtrack},
    $R_{\alpha(i)}$ is disjoint from $B_{\alpha(j)}$, for each $j<i$. Therefore,
    $B'_1,\ldots,B'_n$ is a collection of vertex-disjoint subgraphs. Note also
    that, for each $k\in\{1,\ldots,t\}$ and each $i\in\{2,\ldots,n-1\}$,
    $L_{\alpha(i+1),k}\cap R_{\alpha(i),k}$ is an $(X'_i,Y'_i)$-path intersecting
    $B_{\alpha(i)}$, so $B'_i$ is connected and contains $t$ vertex-disjoint
    $(X'_i,Y'_i)$-paths. Similarly, $B'_1$ and $B'_n$ are each connected and contain
    $t$ vertex-disjoint $(X',Y'_1)$- and $(X'_n,Y')$-paths, respectively.
    Hence, $(B'_1,\ldots,B'_n)$ and $(M'_1,\ldots,M'_{n-1},\emptyset)$ define a
    $(t,0,0,n)$-necklace $N'$, $B'_1$ contains $s$ vertex-disjoint
    $(X,Y'_1)$-paths and $B'_n$ contains $s$ vertex-disjoint $(X'_n,Y)$-paths,
    as desired.
  \end{proof}
\end{lemma}

The following result all but completes the proof of
\autoref{thm:set-to-necklace-preview}.
\begin{lemma}
  \label{lem:grow-theta}
  There exists a function $\fcur:\N^3\rightarrow\N$ such that the following
  holds. Let $t,s,\ell,m,n\in\N$ with $t\geq 1$, $n\geq 2$ and
  $m\geq \fcur(t,s,n)$. Now let $N$ be a long-jump-free $(t,s,\ell,m)$-necklace
  in a graph $G$. If, for each $a,b\in\Z_m$, there exist $(t+s+\ell+1)$
  vertex-disjoint $(V(B_{a}),V(B_{b}))$-paths in $G$, then there is a
  $(t+1,s,\ell,n)$- or $(t,s+1,\ell,n)$-necklace in $G$ supported by $N$.

  \begin{proof}
    Let $t,s,\ell,m,n\in\N$ with $t\geq s$, $t\geq 1$, $n\geq 2$. Let
    \[m_1=\max\{
        \fnum{lem:grow-linear}(t+1,t,n)\,,\,
        \fnum{lem:grow-linear}(s+1,s+1,2)\,,\,
        2s+2\,,\,
        n\}.\]
    Let $a=s+1$, let $b=s+m_1$, and let $\fcur(t,s,n)=2m_1-2$. Now let $N$ be a
    long-jump-free $(t,s,\ell,m)$-necklace in a graph $G$ with
    $m\geq\fcur(t,s,n)$ such that $G$ contains $t+s+\ell+1$ vertex-disjoint
    paths between any two beads; in particular, $G$ contains $t+s+\ell+1$
    vertex-disjoint $(V(B_a),V(B_b))$-paths. By contracting $N$, we may assume
    that $m=2m_1-2$.

    Let $G'=G-Z$ and $N'=(\cB,\cM,\emptyset)$. Note that $G'$ contains $t+s+1$
    vertex-disjoint $(V(B_a),V(B_b))$-paths and that a $(t',s',0)$-necklace
    supported by $N'$ in $G'$ gives a $(t',s',\ell)$-necklace supported by $N$
    in $G$. Therefore, we may assume $\ell=0$.

    Since $N$ is long-jump-free in $G$, there is a separation $(G_1,G_2)$ of $G$
    such that
    \begin{enumerate}[(i)]
    \item $G_1\cap G_2 = B_a\cup B_b$,
    \item $B_a\cup B_{a+1}\cup \cdots \cup B_b$ is a subgraph of $G_1$, and
    \item $B_1\cup B_{2}\cup \cdots \cup B_a$ and
      $B_b\cup B_{b+1}\cup \cdots \cup B_m$ are subgraphs of $G_2$.
    \end{enumerate}

    Let $\cP$ be a collection of $t+s+1$ vertex-disjoint
    $(V(B_a),V(B_b))$-paths; we may assume that each path in $\cP$ meets
    $V(B_a\cup B_b)$ in only its ends. Therefore each path in $\cP$ is contained
    in either $G_1$ or $G_2$; let $\cP_1$ denote the set of paths in $\cP$
    that are contained in $G_1$ and let $\cP_2$ denote the set of paths in $\cP$
    that are contained in $G_2$. Since $|\cP|=t+s+1$, either $|\cP_1|>t$ or
    $|\cP_2|>s$.

    \begin{case-one}
      \quad $|\cP_1|>t$.
    \end{case-one}
    Note that, $(B_a,\ldots,B_b)$ is the bead sequence of a
    $(t,0,0,m_1)$-necklace in $G_1$ that is long-jump-free and has no
    $(1,m_1)$-jump. Moreover, considering the necklace $N$, we see that $G_1$
    contains $t$ vertex-disjoint $(X_a,V(B_{a+t-1}))$-paths and $t$
    vertex-disjoint $(V(B_{b-t+1}),Y_b)$-paths. Since
    $m_1\geq\fnum{lem:grow-linear}(t+1,t,n)$, \autoref{lem:grow-linear} shows that
    $G_1$ contains a $(t+1,0,0,n)$-necklace $N'$ supported by $N$, $B'_1$
    contains $t$ vertex-disjoint $(X_a,Y'_1)$-paths, and $B'_n$ contains
    $(X'_n,Y_b)$-paths.

    By considering the necklace $N$, we see that $G[V(B_1\cup\cdots\cup B_a)]$ contains $s$
    vertex-disjoint $(X_1,X_a)$-paths, each intersecting $G_1$ only at the end
    in $X_a$; adding these paths to $B'_1$ yields a connected subgraph with $s$
    vertex-disjoint $(X_1,Y'_1)$-paths. Similarly, $B'_n$ can be extended to
    contain $s$ vertex-disjoint $(X'_n,Y_m)$-paths, yielding a
    $(t+1,s,0,n)$-necklace supported by $N$.

    \begin{case-two}
      \quad $|\cP_2|>s$.
    \end{case-two}

    If $s=t$, then, because $a+m-b=b-a$, the indices of $N$ can be shifted so
    that $G_1$ and $G_2$ are swapped, reducing to Case 1. We therefore
    assume that $s<t$.

    Note that $(B_b,\ldots,B_m,B_1,\ldots,B_a)$ is the bead sequence of a
    $(0,0,0,m_1)$-necklace in $G_2$ that is long-jump-free and has no
    $(1,m_1)$-jump. Note that $a-s=1$ and $b+s\leq m$, so, considering the
    necklace $N$, $G_2$ contains $s+1$ vertex-disjoint $(X_b,V(B_{b+s}))$-paths
    and $s+1$ vertex-disjoint $(V(B_{a-s}),Y_a)$-paths. Since
    $m_1\geq \fnum{lem:grow-linear}(s+1,s+1,2)$, \autoref{lem:grow-linear} shows
    that $G_2$ contains a $(s+1,0,0,2)$-necklace $N''$ supported by $N$, $B''_1$
    contains $s+1$ vertex-disjoint $(X_b,Y'_1)$-paths, and $B''_2$ contains
    $s+1$ vertex-disjoint $(X'_2,Y_a)$-paths. Then
    $(B''_2,B_{a+1},\ldots,B_{b-1},B''_1)$ is the bead sequence of a
    $(t,s+1,0,m_1)$-necklace supported by $N$; since $m_1\geq n$,
    a $(t,s+1,0,n)$-necklace can be obtained by contraction.
  \end{proof}
\end{lemma}

We can now prove \autoref{thm:set-to-necklace-preview} which we restate below.

\begin{theorem}
  \label{thm:necklace-main}
  There exists a function $\fcur:\N^2\rightarrow\N$ such that the following
  holds. Let $\theta,m,n\in\N$ with $n\ge 2$ and $m\ge \fcur(\theta,n)$. If $G$
  is a graph containing a $\theta$-connected set $U$ of size at least $m$, then
  $G$ contains a $(\theta;n)$-necklace supported by $U$.

  \begin{proof}
    Fix $\theta,n\in\N$, and define
    \[\fcur(\theta,n)=\fnum{lem:init-necklace}(\theta,n,g(1,0,0)),\]
    where, for $t,s,\ell,k\in\N$,
    \[\begin{aligned}
        g(t,s,\ell)&=\begin{cases}
          n & \mbox{ if } t+s+\ell\geq\theta, \\
          \fnum{lem:long-jump-free}(s,\fnum{lem:grow-theta}(t,s,h(t+s+\ell+1)) & \mbox{
            otherwise }
        \end{cases} \\
        h(k)&=\max\{g(t',s',\ell')\,:\,t',s',\ell'\in\N,\,t'\geq 1,\,t'+s'+\ell'=k\}.
      \end{aligned} \]
    Note that this mutual recurrence is strictly decreasing in
    $\theta-t-s-\ell$, and hence well defined. Let $G$ be a graph with a
    $\theta$-connected set $U$ of size at least $\fcur(\theta,n)$.

    By \autoref{lem:init-necklace}, $G$ contains either a
    $(0,0,\theta,n)$-necklace supported by $U$ or a $(1,0,0,g(1,0,0))$-necklace
    such that each bead contains at least $\theta$ vertices in $U$. In the
    former case, $G$ contains a $(\theta;n)$-necklace as required. Suppose now
    that $G$ contains a $(1,0,0,g(1,0,0))$-necklace $N^0$ supported by $U$.
  
    Let $t,s,\ell\in\N$ such that $t\geq 1$, $G$ contains a
    $(t,s,\ell,g(t,s,\ell))$-necklace $N$ supported by $N^0$ and, subject to
    this, $t+s+\ell$ is maximized. We may assume that $t+s+\ell<\theta$.

    Let $m''=h(t+s+\ell+1)$, let $m'=\fnum{lem:grow-theta}(t,s,m'')$, and let
    $m=\fnum{lem:long-jump-free}(s,m')$. Note that $m'=g(t,s,\ell)$. Note also
    that $m'\geq m''$ because the hypothesis of \autoref{lem:grow-theta} is
    satisfied for some graph $H$ and and some $(t,s,\ell,m')$-necklace in $H$,
    so $H$ contains a $(t+1,s,\ell,m'')$- or $(t,s+1,\ell,m'')$-necklace
    supported by a $(t,s,\ell,m')$-necklace.

    By \autoref{lem:long-jump-free}, $G$ contains a $(t',s',\ell',m')$-necklace
    $N'$ supported by $N$ such that $t'\geq 1$, $t'+s'+\ell'\geq t+s+\ell$ and,
    if equality holds, then $N'$ is long-jump-free. If $t'+s'+\ell'>t+s+\ell$,
    then, by deleting some edges, we may assume $t'+s'+\ell'=t+s+\ell+1$;
    because $m'\geq m''\geq g(t',s',\ell')$, $G$ contains a
    $(t',s',\ell',g(t',s',\ell'))$-necklace supported by $U$, contradicting our
    choice of $t,s,\ell$.

    Suppose now that $t'+s'+\ell'=t+s+\ell$, and $N'$ is long-jump-free. Note
    that, for $a,b\in\Z_{m'}$, $V(B'_a)$ and $V(B'_b)$ each contain $\theta$
    vertices in $U$, so $G$ contains $\theta$ vertex-disjoint
    $(V(B'_a),V(B'_b))$-paths. By \autoref{lem:grow-theta}, $G$ contains a
    $(t'+1,s',\ell',m'')$- or $(t',s'+1,\ell',m'')$-necklace. However,
    $m''\geq \max\{g(t'+1,s',\ell'),g(t',s'+1,\ell')\}$, contradicting our
    choice of $t,s,\ell$.
  \end{proof}
\end{theorem}

\section{Nearly-balanced necklaces}
\label{sec:balancing}

A $(t,s,\ell,n)$-necklace $N$ is \emph{balanced} if $s=t$; the necklace $N$ is
\emph{nearly-balanced} if either $N$ is balanced or $s= t-1$ and $\ell=0$. We
show in this section that a $(\theta;m)$-necklace gives rise to a
nearly-balanced $(\theta;n)$-necklace in the same graph; the same is not generally
true of balanced necklaces\textemdash if $\theta$ is odd, then any balanced
$(\theta;n)$-necklace must have at least one hub, which must have degree at
least $n$, so graphs with maximum degree less than $n$ do not have any balanced
$(\theta;n)$-necklace. In the next section we show that a $(\theta;m)$-necklace
does give rise to a balanced $(\theta;n)$-necklace in a minor of the original graph.

The following lemma shows that certain jumps can be used to increase $s$, the
size of the matching between the first and last bead. Such jumps will arise in
two different ways, depending on whether our necklace has a hub
(\autoref{lem:remove-hub}) or not (\autoref{lem:extract-path}).

\begin{lemma}
  \label{lem:jump-increase-s}
  For $t,s,\ell,n,m\in\N$ with $t>s$, $n\geq 2$ and $m\geq n+2s$, if $N$ is a
  $(t,s,\ell,m)$-necklace with an $(s+1,s+n)$-jump in a graph $G$, then $G$
  contains a $(t,s+1,\ell,n)$-necklace supported by $N$.
  \begin{proof}
    By contracting $N$, we may assume $m=n+2s$. By growing the beads of $N$, we
    may assume that $G$ contains an $(s+1,s+n)$-jump consisting of exactly one
    edge, $e$, incident with a vertex $x\in V(B_{s+1})$ and a vertex in
    $y \in V(B_{s+n})$. Let
    $$N' = N\circ (s+1,s+2,\ldots,s+n-1).$$
    Now we construct $N''$ from $N'$ by adding $e$ to $M_n$. Let
    $\cX''=(X''_1,\ldots,X''_{n})$ and $\cY''=(Y''_1,\ldots, Y''_{n})$ be the
    left and right attachment sequences for $N''$.

    We claim that $N''$ is a $(t,s+1,\ell,n)$-necklace in $G$. By The only
    non-trivial details that need to be verified are:
    \begin{enumerate}[(1)]
    \item there are $s+1$ vertex-disjoint $(X''_1,Y''_{1})$-paths in $B''_1$,
      and
    \item there are $s+1$ vertex-disjoint $(X''_{n},Y''_{n})$-paths in
      $B''_{n}$.
    \end{enumerate}
    Note that $X''_1=X_1\cup\{x\}$, $Y''_1=Y_{s+1}$, $X''_n=X_{s+n}$ and
    $Y''_n=Y_m\cup\{y\}$. Note also that $m-(s+n)=s$, so reversing $N$ and
    swapping $x$ and $y$ reduces $(2)$ to $(1)$. Therefore, it suffices to prove
    $(1)$.

    Considering the necklace $N$, $B''_1$ contains a collection $\cP$ of $t$
    vertex-disjoint $(Y_1,Y_{s+1})$-paths such that, for each $P\in\cP$ and each
    $i\in\{2,\ldots,s+1\}$, $P\cap B_i$ is an $(X_i,Y_i)$-path. Choose
    $P_0\in\cP$ such that $B_{s+1}$ contains an $(x,V(B_{s+1}\cap P_0))$-path
    $P'$ disjoint from each path in $\cP-P_0$. Let $P'_0$ be an
    $(x,Y_{s+1})$-path in $P'\cup P_0\cap B_{s+1}$. Because $t>s$, there exist
    distinct paths $P_1,\ldots,P_s\in \cP-\{P_0\}$.

    Considering the necklace $N$ again, $B''_1$ contains a collection
    $(Q_1,\ldots,Q_s)$ of vertex-disjoint $(X_1,X_s)$-paths, each intersecting
    each of $B_1,\ldots,B_s$, and each disjoint from $B_{s+1}$. Let
    $H=\bigcup_{i=1}^{s}(B_i\cup P_i\cup Q_i)$.
    By \autoref{lem:connectivity}, $H$ contains a collection
    $(P'_1,\ldots,P'_s)$ of vertex disjoint $(X_1,Y_{s+1})$-paths. Moreover,
    $V(H\cap P'_0)=\emptyset$, so $(P'_0,\ldots,P'_s)$ is a collection of
    $s+1$ vertex-disjoint $(X_1\cup\{x\},Y_{s+1})$-paths. Hence $(1)$ holds.
  \end{proof}
\end{lemma}

Note that a hub can be turned into a jump between any two beads, so
\autoref{lem:jump-increase-s} immediately implies that a hub can be used to
increase $s$, the size of the matching between the first and last bead.
\begin{lemma}
  \label{lem:remove-hub}
  Let $t,s,\ell,n,m\in\N$ with $t>s$, $\ell\geq 1$, $n\geq 2$ and $m\ge n+2s$.
  If $N$ is a $(t,s,\ell,m)$-necklace in a graph $G$, then there is a
  $(t,s+1,\ell-1,n)$-necklace in $G$ that is supported by $N$.
\end{lemma}

The following result handles unbalanced necklaces without hubs; in this case we
remove from the necklace one of the $t$ paths through the beads
$B_1,\ldots,B_n$. This path gives a jump between any two beads, allowing
\autoref{lem:jump-increase-s} to be applied. This is also used in the next
section, where the path is contracted to form a new hub.
\begin{lemma}
  \label{lem:extract-path}
  There exists a function $\fcur:\N^2\rightarrow \N$ such that the following
  holds. Let $t,s,m,n\in\N$ with $t>s$, $t+s\geq 2$, $n\geq 2$ and
  $m\geq \fcur(t,n)$. If $N$ is a $(t,s,0,m)$-necklace supported by a set $U$ of
  vertices in a graph $G$, then $G$ contains a $(t-1,s,0,n)$-necklace $N'$
  supported by $U$, and a path $P$ such that $V(P)\cap V(N')=\emptyset$ and, for
  $i\in\{1,\ldots,n\}$, $G$ contains an edge between $V(B'_i)$ and $V(P)$.

  \begin{proof}
    Let $t,s,n\in\N$ with $t> s$, $t+s\geq 2$, and $n\geq 2$. Let
    $\fcur(t,n)=t(n-1)+1$, and let $m\geq \fcur(t,n)$. Let $N$ be a
    $(t,s,0,m)$-necklace supported by a vertex set $U$ in a graph $G$.

    Let $P_1,\ldots,P_t$ be a collection of vertex-disjoint $(X_1,Y_m)$-paths in
    $G[N]\delete M_m$. Note that, for $i\in\Z_m$ and $k\in\{1,\ldots,t\}$,
    $P_k\cap B_i$ is an $(X_i,Y_i)$-path. If $B_i-V(P_k)$ contains a connected
    component $B'$ such that, for $j\in\{1,\ldots,t\}-\{k\}$, $V(P_j\cap B_i)\subseteq V(B')$
    and $U\cap V(B')\neq\emptyset$, then we
    say that $P_k$ is \emph{removable at $i$}.

    \begin{claim}
      For each $i\in\Z_m$, there exists $k\in\{1,\ldots,t\}$ such that $P_k$ is
      removable at $i$.
      \begin{subproof}
        Fix $i\in\Z_m$ and choose $u\in U\cap V(B_i)$. For $k\in\{1,\ldots,t\}$,
        let $(H^k_1,H^k_2)$ be a separation in $B_i$ such that
        $H_1^k\cap H_2^k=P_k\cap B_i$, $u\in V(H^k_1)$, and $H^k_1-V(P_j)$ is
        connected (but possibly empty). Note that $2t>t+s\geq 2$, so $t>1$.
        Therefore, $u\notin V(P_k)$ for some $k\in\{1,\ldots,t\}$. Choose
        $k\in\{1,\ldots,t\}$ such that $u\notin V(P_k)$ and $V(H^k_1)$ is
        maximal. Suppose for contradiction that $P_k$ is not removable at $i$,
        so there exists $j\in\{1,\ldots,t\}-\{k\}$ such that
        $V(P_j\cap B_i)\not\subseteq V(H^k_1)$. Since $P_j\cap B_i$ is connected
        and disjoint from $V(H^k_1\cap H^k_2)$,
        $V(P_j\cap B_i)\subseteq V(H^k_2)-V(H^k_1)$.
        Therefore, $H^k_1\cup (P_k\cap B_i)$ is a
        connected subgraph of $B_i-V(P_j)$ containing $u$, so $u\notin V(P_j)$
        and $V(H^k_1)\subseteq V(H^k_1\cup (P_k\cap B_i))\subsetneq V(H^j_1)$, contradicting the choice of $k$.
      \end{subproof}
    \end{claim}

    Therefore, there exist $k\in\{1,\ldots,t\}$ and a set
    $I\subseteq \{1,\ldots,m\}$ with $|I|\ge n$ such that for each $i\in I$,
    $P_k$ is removable at $i$; without loss of generality, $k=t$. Let
    $(a_1,\ldots, a_{n})$ be an increasing sequence of integers in $I$. Let $P$
    be the subpath of $P_t$ between $X_{a_1+1}$ and $Y_{a_n-1}$. For
    $i\in\{1,\ldots,n\}$, let $\{e_i\}=E(M_{a_i}\cap P_t)$. Let
    $N'= N\circ (a_1,\ldots, a_{n-1})$. Let $B''_1=B'_1$, let $B''_n=B'_n$. For
    $i\in\{2,\ldots,n-1\}$, $P_t$ is removable at $i$ in $N'$; let $B''_i$ be
    the component of $B'_i-V(P)$ containing $P_j\cap B'_i$, for each
    $j\in\{1,\ldots,t-1\}$, and containing a vertex in $U$.
    For $i\in\{1,\ldots,n-1\}$, let $M''_i=M'_i-\{e_i\}$.
    Then $N''=((B''_1,\ldots,B''_n),(M''_1,\ldots,M''_n),\emptyset)$ is a
    $(t-1,s,0,n)$-necklace supported by $U$.
    Finally, $e_1$ is an edge between $V(P)$ and $V(B''_1)$, $e_{n-1}$ is an
    edge between $V(P)$ and $V(B''_n)$ and for $i\in\{2,\ldots,n-1\}$, $B'_i$
    contains an edge between $V(P)$ and $V(B''_i)$.
  \end{proof}
\end{lemma}

The path extracted by \autoref{lem:extract-path} gives a jump between any two
beads. By \autoref{lem:jump-increase-s}, such a jump can be used to increase $s$,
implying the following:

\begin{lemma}
  \label{lem:flip-path}
  Let $t,s,m,n,\in\N$ with $t\geq s+2$, $n\geq 2$ and
  $m\geq \fnum{lem:extract-path}(t,n+2s)$. If $N$ is a $(t,s,0,m)$-necklace
  supported by a set $U$ in a graph $G$, then there is a
  $(t-1,s+1,0,n)$-necklace in $G$ supported by $U$.
\end{lemma}

Lemmas~\ref{lem:remove-hub} and~\ref{lem:flip-path}, together with a routine
inductive argument on $(t-s)$ (which we omit), show that a $(t,s,\ell,m)$
necklace gives rise to a nearly-balanced $(t+s+\ell;n)$-necklace. Combined with
\autoref{thm:necklace-main}, this proves the following, showing that large
$\theta$-connected sets give rise to nearly-balanced $(\theta;n)$-necklaces.

\begin{theorem}
  \label{thm:nearly-balanced-necklace}
  There exists a function $\fcur:\N^2\rightarrow \N$ such that the following
  holds.
  Let $\theta,n\in\N$ with $n\geq 2$.
  If $U$ is a $\theta$-connected vertex set in a graph $G$ and $|U|\geq\fcur(\theta,n)$,
  then $G$ contains a nearly-balanced $(\theta,n)$-necklace supported by $U$.


\end{theorem}

In some applications, this result may be more useful than
\autoref{thm:main-local} because the necklace is a subgraph, rather than a
minor.

\section{Generalized wheel-minors from necklaces}
\label{sec:necklace-to-wheel}

The penultimate step of the proof is to find a balanced $(\theta;n)$-necklace in
a minor of $G$. The final step, getting a generalized wheel from a balanced
necklace, is straightforward.

A ``retract'' of $G$ is a special kind of contraction-minor in which we take
greater care of vertex labels. If $e=xy$ is an edge of a graph $G$, then the
{\em retraction of $e$ to $x$ in $G$} is the graph obtained by contracting $e$
to the vertex $x$; so there are two ways of retracting a non-loop edge. A {\em
  retract} of $G$ is any graph obtained from $G$ by a, possibly empty, sequence
of edge retractions. Therefore, if $H$ is a retract of $G$, then
$V(H)\subseteq V(G)$. A graph $H$ is a {\em retraction-minor} of a graph $G$ if
$H$ is a subgraph of a retract of $G$.

A $(t,s,\ell,n)$-necklace $N$ is {\em refined} if each bead is a tree with
$\max\{1,t\}$ vertices. Given a refined $(t,s,\ell,n)$-necklace $N$ we will
assume, for each $i\in \Z_n$, that $V(B_i) = \{x_{i,1},\ldots, x_{i,t'}\}$ and,
for $i\in\{1,\ldots,n-1\}$, that
$M_i = \{x_{i,1}x_{i+1,1},\ldots, x_{i,t'}x_{i+1,t'}\}$, where $t'=\max\{1,t\}$.

We omit the routine proof of the following result.
\begin{lemma}
  \label{lem:refined-retract}
  If $N$ is a $(t,s,\ell,n)$-necklace supported by a set $U$ of vertices in a
  graph $G$, then there is a retract $H$ of $G$ and a refined
  $(t,s,\ell,n)$-necklace supported by $U$ in $H$.
\end{lemma}

The following result gives a balanced necklace in a minor.
\begin{lemma}\label{lem:balanced-necklace}
  There exists a function $\fcur:\N^2\rightarrow \N$ such that the following
  holds. Let $\theta,n\in\N$ with $n\geq 2$. If $G$ is a graph with a
  $\theta$-connected set $U$ of size at least $\fcur(\theta,n)$, then there is a
  retract $H$ of $G$ containing a balanced $(\theta;n)$-necklace supported by
  $U$.
  \begin{proof}
    Fix $\theta,n\in\N$ and define
    \[\begin{aligned}
        g(\theta,n) &= \max\{n, \fnum{lem:extract-path}(\lfloor  \theta/2 \rfloor,n)\} \\
        \fcur(\theta,n) &= \fnum{thm:nearly-balanced-necklace}(\theta,g(\theta,n)).
      \end{aligned} \] Let $G$ be a graph with a $\theta$-connected set $U$ of
    size at least $\fcur(\theta,n)$. By \autoref{thm:nearly-balanced-necklace},
    $G$ contains a nearly-balanced $(t,s,\ell,g(\theta,n))$-necklace $N$
    supported by $U$ with $t+s+\ell=\theta$. If $N$ is balanced, then $G$ itself
    contains a balanced $(\theta;n)$-necklace. Suppose $N$ is not balanced, so
    $N$ is a $(s+1,s,0,g(\theta,n))$-necklace.
    By
    \autoref{lem:extract-path}, $G$ contains $(s,s,0,n)$-necklace $N'$ supported
    by $U$ and a path $P$ such that $V(P)\cap V(N')=\emptyset$ and each bead in
    $N'$ is adjacent to a vertex in $P$. Let $H$ be the retract-minor of $G$
    obtained contracting $E(P)$ to some vertex $z\in V(P)$. Then $N'$ is a
    necklace supported by $U$ in $H$, and $z$ is adjacent to each bead of $N'$,
    so $H$ contains an $(s,s,1,n)$-necklace supported by $U$.
  \end{proof}
\end{lemma}

The main result of the paper is:
\begin{theorem}
  \label{thm:main-local}
  There exists a function $\fcur:\N^2\rightarrow\N$ such that the following
  holds. Let $\theta,n\in\N$ with $\theta\geq 2$ and $n\geq 3$. If $G$ is a
  graph with a $\theta$-connected set $U$ of size at least $\fcur(\theta,n)$,
  then $G$ has a retraction minor $H$ such that either
  \begin{enumerate}[(i)]
  \item $H$ is isomorphic to $K_{\theta,n}$ with an independent set of size $n$
    contained in $U$, or
  \item $H$ is a $(\theta; n)$-wheel with a rim-transversal contained in $U$.
  \end{enumerate}

  \begin{proof}
    Fix $\theta,n\in\N$. Define
    $m=(n-1) \left\lfloor \theta/2 \right\rfloor^{\lfloor \theta/2 - 2\rfloor}+1$,
    and define
    $\fcur(\theta,n)=\fnum{lem:balanced-necklace}\left(\theta, m\right)$. Let
    $U$ be a $\theta$-connected vertex set in a graph $G$ with
    $|U|\geq\fcur(\theta,n)$.

    By Lemmas~\ref{lem:balanced-necklace} and~\ref{lem:refined-retract}, there
    is a retract $H_1$ of $G$ and a refined $(t,t,\ell,m)$-necklace $N$ of $H_1$
    that is supported by $U$ with $2t+\ell = \theta$. If $t=0$, then each bead
    of $N$ contains exactly one vertex, which is in $U$, so $H_1[N]$ is a
    complete bipartite graph satisfying $(i)$. We now assume that $t>0$.

    For each $i\in \Z_n$, we let $T_i$ be an isomorphic copy of $B_i$ obtained
    by ``renaming'' vertex $x_{i,j}$ as $x_j$, for $j\in\{1,\ldots,t\}$. Since
    there are $t^{t-2}$ trees on the vertex set $\{x_1,\ldots,x_t\}$, and
    $(n-1)t^{t-2}<m$ there is a $n$-element set $I\subseteq \{1,\ldots, m\}$
    such that the trees $(T_i\, : \, i\in I)$ are all equal. It is now routine
    to obtain the $(\theta,n)$-wheel required in $(ii)$.
  \end{proof}
\end{theorem}

\section{The Grid Theorem}
\label{sec:gridtheorem}

We conclude by proving Robertson and Seymour's Grid Theorem using
\autoref{thm:main-local}.

\begin{theorem}[Grid Theorem]
  \label{thm:grid-on-set}
  For each natural number $n$ there exists a natural number $m$ such that, if
  $G$ is a graph containing an $m$-connected set of size at least $m$, then $G$
  contains an $n\times n$-grid-minor.

  \begin{proof}
    Define $\theta = 2\cdot n^{2n}+n^2$ and
    $m=\fnum{thm:main-local}(\theta,n^2)$. Let $G$ be a graph with an
    $m$-element $m$-connected set. Since the $n\times n$-grid is a subgraph of
    $K_{n^2,n^2}$, we may assume that $G$ has no $K_{n^2,n^2}$-minor. Then,
    by~\autoref{thm:main-local}, $G$ has a $(t,\ell,n^2)$-wheel $H$ as a minor,
    where $2t+\ell=\theta$. Suppose that the wheel $H$ is defined by
    $(T,Z,\pi,\psi)$. Since $G$ has no $K_{n^2,n^2}$-minor, $\ell<n^2$ and,
    hence, $t > n^{2n}$.

    We may assume that $T$ does not contain an $n$-vertex path, since
    otherwise $H$ would contain a subgraph isomorphic to the $n\times n$-grid.
    Since $|V(T)|=t> (n^2)^n$, there is a vertex of $v$ degree at least $n^2$ in
    $T$. Now, contracting the edges $u_iu_{i+1}$ for each neighbour $u$ of $v$
    in $T$ and each $i\in\{1,\ldots,n^2-1\}$, we obtain a $K_{n^2,n^2}$-minor;
    this contradiction completes the proof.
  \end{proof}
\end{theorem}

\bibliographystyle{plain}
\bibliography{library}

\end{document}